\documentclass[12pt]{article}

\usepackage{amssymb,latexsym}

\newtheorem{theorem}{Theorem}[section]
\newtheorem{proposition}[theorem]{Proposition}
\newtheorem{corollary}[theorem]{Corollary}
\newtheorem{definition}[theorem]{Definition}

\newtheorem{lemma}[theorem]{Lemma}

\newtheorem{conjecture}[theorem]{Conjecture}

\def\cM{{\cal M}}
\def\cN{{\cal N}}
\def\cO{{\cal O}}

\newenvironment{proof}{\medskip\noindent{\it Proof.\ }}{\mbox{$\Box$}\medskip}

\def\cplus{\ {{\ } \atop {+}}\ }
\def\cminus{\ {{\ } \atop {-}}\ }

\begin{document}

\title{Involutions Restricted by 3412, Continued Fractions, and Chebyshev Polynomials}

\author{Eric S. Egge \\
Department of Mathematics \\
Gettysburg College\\
Gettysburg, PA  17325  USA \\[4pt]
eggee@member.ams.org \\
\\
Toufik Mansour\\
Department of Mathematics, Haifa University\\
31905 Haifa, Israel \\[4pt]
toufik@math.haifa.ac.il}
\maketitle

\begin{abstract}
We study generating functions for the number of involutions, even
involutions, and odd involutions in $S_n$ subject to two
restrictions. One restriction is that the involution avoid $3412$
or contain 3412 exactly once. The other restriction is that the
involution avoid another pattern $\tau$ or contain $\tau$ exactly
once. In many cases we express these generating functions in terms
of Chebyshev polynomials of the second kind.

\medskip

\end{abstract}

\section{Introduction and Notation}

Let $S_n$ denote the set of permutations of $\{1, \ldots, n\}$, written in one-line notation, and suppose $\pi \in S_n$.
We write $|\pi|$ to denote the length of $\pi$, and for all $i$, $1 \le i \le n$, we write $\pi(i)$ to denote the $i$th element of $\pi$.
We say $\pi$ is an {\it involution} whenever $\pi(\pi(i)) = i$ for all $i$, $1 \le i \le n$, and we write $I_n$ to denote the set of involutions in $S_n$.
Now suppose $\pi \in S_n$ and $\sigma \in S_k$.
We say a subsequence of $\pi$ has {\it type} $\sigma$ whenever it has all of the same pairwise comparisons as $\sigma$.
For example, the subsequence 2869 of the permutation 214538769 has type 1324.
We say $\pi$ {\it avoids} $\sigma$ whenever $\pi$ contains no subsequence of type $\sigma$.
For example, the permutation 214538769 avoids 312 and 2413, but it has 2586 as a subsequence so it does not avoid 1243.
In this context $\sigma$ is sometimes called a {\it pattern} or {\it forbidden subsequence} and $\pi$ is sometimes called a {\it restricted permutation} or {\it pattern-avoiding permutation}.
In this paper we will be interested in involutions which avoid several patterns, so for any set $R$ of permutations we write $I_n(R)$ to denote the set of involutions in $S_n$ which avoid every pattern in $R$ and we write $I(R)$ to denote the set of all involutions which avoid every pattern in $R$.

In \cite{Egge} Egge connected generating functions for various subsets of $I(3412)$ with continued fractions and Chebyshev polynomials of the second kind, and gave a recursive formula for computing them.
For instance, he showed \cite[Theorem 3.3]{Egge} that
\begin{equation}
\label{eqn:introcf}
\sum_{\pi \in I(3412)} \prod_{k \ge 1} x_k^{\tau_k(\pi)} = \frac{1}{\displaystyle 1 - x_1 - \frac{x_1^2 x_2}{\displaystyle
1 - x_1 x_2^2 x_3 - \frac{x_1^2 x_2^5 x_3^4 x_4}{\displaystyle
1 - x_1 x_2^4 x_3^6 x_4^4 x_5 - \cdots}}}.
\end{equation}
Here the $n$th numerator is $\prod\limits_{i=1}^{2n} x_i^{{{2n-2}\choose{i-1}} + {{2n-1}\choose{i-1}}}$ and the $n$th denominator is $1 - \prod\limits_{i=1}^{2n+1} x_i^{{{2n}\choose{i-1}}}$.
Moreover, we have used the notation we will adopt throughout the paper:  we write $[k]$ to denote the permutation $k \ldots 21$ and we write $\tau_k(\pi)$ to denote the number of subsequences of type $[k]$ in $\pi$.
Egge also showed \cite[Theorem 6.1]{Egge} that for all $k \ge 1$,
\begin{equation}
\sum_{\sigma \in I(3412, [2k])} x^{|\sigma|} = \frac{U_{k-1}\left(\frac{1-x}{2x}\right)}{x U_k\left(\frac{1-x}{2x}\right)}
\end{equation}
and
\begin{equation}
\sum_{\sigma \in I(3412, [2k-1])} x^{|\sigma|} = \frac{ U_{k-1}\left(\frac{1-x}{2x}\right) + U_{k-2}\left(\frac{1-x}{2x}\right)}{x\left( U_k\left(\frac{1-x}{2x}\right) + U_{k-1}\left(\frac{1-x}{2x}\right)\right)}.
\end{equation}
Here $U_n(x)$ is the $n$th Chebyshev polynomial of the second kind, which may be defined by ${\displaystyle U_n(\cos t) = \frac{\sin((n+1) t)}{\sin t}}$.
Egge proved these results using the following recursive formula \cite[Corollary 5.6]{Egge} for the generating function ${\displaystyle F^+_\pi(x) = \sum_{\sigma \in I(3412,\pi)} x^{|\sigma|}}$, which makes it possible to compute $F^+_\pi(x)$ for any permutation $\pi$:
\begin{equation}
\label{eqn:introrecurrence}
F^+_\pi(x) = 1 + x F^+_\beta(x) + x^2 \sum_{i=1}^k \left( F^+_{\overline{\alpha_1 \oplus \cdots \oplus \alpha_i}}(x) - F^+_{\overline{\alpha_1 \oplus \cdots \oplus \alpha_{i-1}}}(x)\right) F^+_{\alpha_i \oplus \cdots \oplus \alpha_k}(x).
\end{equation}
Here the various subscripts of $F^+$ on the right are the types of certain subsequences of $\pi$.
For other results concerning pattern-avoiding permutations, continued fractions, and Chebyshev polynomials, see \cite{EggeMansourSchroder,MansourVainshtein5,ARe} and the references therein.

In this paper we refine Egge's results by studying generating functions for even and odd involutions in $I(3412)$.
For any permutation $\pi$, we write $sign(\pi)$ to denote the sign of $\pi$, which is 1 if $\pi$ is an even permutation and $-1$ if $\pi$ is an odd permutation.
Using Egge's techniques, we prove signed analogues of (\ref{eqn:introcf})--(\ref{eqn:introrecurrence}), from which one can obtain analogues for even and odd permutations.
For instance, we show that
\begin{displaymath}
\sum_{\pi \in I(3412)} (-1)^{sign(\pi)} \prod_{k \ge 1} x_k^{\tau_k(\pi)} = \frac{1}{\displaystyle 1 - x_1 + \frac{x_1^2 x_2}{\displaystyle
1 - x_1 x_2^2 x_3 + \frac{x_1^2 x_2^5 x_3^4 x_4}{\displaystyle
1 - x_1 x_2^4 x_3^6 x_4^4 x_5 + \cdots}}}.
\end{displaymath}
Here the $n$th numerator is $\prod\limits_{i=1}^{2n} x_i^{{{2n-2}\choose{i-1}} + {{2n-1}\choose{i-1}}}$ and the $n$th denominator is $1 - \prod\limits_{i=1}^{2n+1} x_i^{{{2n}\choose{i-1}}}$.
We also show that for all $k \ge 1$,
\begin{displaymath}
\sum_{\sigma \in I(3412, [2k])} (-1)^{sign(\pi)} x^{|\sigma|} = \frac{U_{k-1}\left(\frac{1-x}{2ix}\right)}{ix U_k\left(\frac{1-x}{2ix}\right)}
\end{displaymath}
and
\begin{displaymath}
\sum_{\sigma \in I(3412, [2k-1])} (-1)^{sign(\pi)} x^{|\sigma|} = \frac{ U_{k-1}\left(\frac{1-x}{2ix}\right) - i U_{k-2}\left(\frac{1-x}{2ix}\right)}{ix\left( U_k\left(\frac{1-x}{2ix}\right) - i U_{k-1}\left(\frac{1-x}{2ix}\right)\right)}.
\end{displaymath}
These results follow from our analogue of (\ref{eqn:introrecurrence}), which is the following recursive formula for the signed generating function ${\displaystyle F^-_\pi(x) = \sum_{\sigma \in I(3412, \pi)} (-1)^{sign(\sigma)} x^{|\sigma|}}$:
\begin{displaymath}
F^-_\pi(x) = 1 + x F^-_\beta(x) - x^2 \sum_{i=1}^k \left( F^-_{\overline{\alpha_1 \oplus \cdots \oplus \alpha_i}}(x) - F^-_{\overline{\alpha_1 \oplus \cdots \oplus \alpha_{i-1}}}(x)\right) F^-_{\alpha_i \oplus \cdots \oplus \alpha_k}(x).
\end{displaymath}
This result makes it possible to compute $F^-_\pi(x)$ for any permutation $\pi$.

In addition to signed generating functions for various subsets of $I(3412)$, we also study generating functions for involutions which contain a given pattern exactly once and avoid another pattern.
For any patterns $\sigma_1$ and $\sigma_2$ and any $n \ge 0$, we write $I_n(\sigma_1; \sigma_2)$ (resp. $I(\sigma_1; \sigma_2)$) to denote the set of involutions of length $n$ (resp. of any length) which avoid $\sigma_1$ and contain exactly one subsequence of type $\sigma_2$.
We first use Egge's description of the recursive structure of $I(3412)$ to obtain recurrence relations for the generating functions ${\displaystyle G^+_\pi(x) = \sum_{\sigma \in I(3412; \pi)} x^{|\sigma|}}$ and ${\displaystyle G^-_\pi(x) = \sum_{\sigma \in I(3412; \pi)} (-1)^{sign(\sigma)} x^{|\sigma|}}$ for certain $\pi$.
Using these recurrence relations we find $G^+_\pi(x)$ and $G^-_\pi(x)$ for various $\pi$ in terms of Chebyshev polynomials of the second kind.
For example, we show that for all $k \ge 3$,
\begin{displaymath}
G^+_{k\ldots 4213}(x) = \frac{1-x-x^2}{x\left( U_{k-1}\left(\frac{1-x}{2x}\right) + U_{k-2}\left(\frac{1-x}{2x}\right)\right)^2}
\end{displaymath}
and
\begin{displaymath}
G^-_{k\ldots 4213}(x) = \frac{x (1 - x + x^2)}{\left(U_{k-2}\left(\frac{1-x}{2ix}\right) - ix U_{k-3}\left(\frac{1-x}{2ix}\right)\right)^2},
\end{displaymath}
and we find similar results when $213$ is replaced with any permutation of length three.

We then turn our attention to the set $I(;3412)$ of involutions which contain exactly one subsequence of type 3412.
We first describe the recursive structure of $I(;3412)$, which we use to find the signed and unsigned generating functions for $I(;3412)$.
We then use this structure to obtain recurrence relations for the generating functions ${\displaystyle P^+_\pi(x) = \sum_{\sigma \in I(\pi;3412)} x^{|\sigma|}}$ and ${\displaystyle P^-_\pi(x) = \sum_{\sigma \in I(\pi;3412)} (-1)^{sign(\sigma)} x^{|\sigma|}}$.
These recurrence relations enable us to find $P^+_\pi(x)$ and $P^-_\pi(x)$ for various $\pi$ in terms of Chebyshev polynomials of the second kind.
For example, we show that for all $k \ge 1$ we have
\begin{displaymath}
P^+_{[2k]}(x) = \frac{\sum_{j=0}^{k-2}(1-x^{2j+2}) U_j^2\left(\frac{1-x}{2x}\right)}{(1-x) U_k^2\left(\frac{1-x}{2x}\right)}
\end{displaymath}
and
\begin{displaymath}
P^-_{[2k]}(x) = \frac{\sum_{j=0}^{k-1} v_{2j}(x) U_j^2\left(\frac{1-x}{2ix}\right)}{U_k^2\left(\frac{1-x}{2ix}\right)},
\end{displaymath}
where $v_k(x) = \sum_{j=0}^k (-1)^{{{j}\choose{2}}} x^j$.

We conclude the paper with a short list of directions for future research.

\section{The Recursive Structure of $I(3412)$}

In this section we recall the natural recursive structure of $I(3412)$, which was first observed by Guibert \cite[Remark 4.28]{G}.
We begin with notation for three ways of combining permutations.

\begin{definition}
Suppose $\pi \in S_m$ and $\sigma \in S_n$.
We write $\pi \oplus \sigma$ to denote the permutation in $S_{m+n}$ given by
$$(\pi \oplus \sigma)(i) =
\cases{
\pi(i) & if $1 \le i \le m$, \cr
\sigma(i-m)+m & if $m+1 \le i \le m+n$.
}
$$
We refer to $\pi \oplus \sigma$ as the {\em direct sum} of $\pi$ and $\sigma$.
\end{definition}

\begin{definition}
Suppose $\pi \in S_m$ and $\sigma \in S_n$.
We write $\pi \ominus \sigma$ to denote the permutation in $S_{m+n}$ given by
$$(\pi \ominus \sigma)(i) =
\cases{
\pi(i) + n & if $1 \le i \le m$, \cr
\sigma(i-n) & if $m+1 \le i \le m+n$.}
$$
We refer to $\pi \ominus \sigma$ as the {\em skew sum} of $\pi$ and $\sigma$.
\end{definition}

\begin{definition}
Suppose $\pi \in S_m$ and $\sigma \in S_n$.
We write $\pi * \sigma$ to denote the permutation in $S_{m+n+2}$ given by
$$\pi * \sigma = (1 \ominus \pi \ominus 1) \oplus \sigma.$$
\end{definition}

We now describe the recursive structure of $I(3412)$.

\begin{proposition}
\label{prop:permmap}
(\cite[Proposition 2.8]{Egge})
\renewcommand\labelenumi{{\upshape (\roman{enumi}) }}
\begin{enumerate}
\item
For all $n \ge 1$, the map
$$
\begin{array}{ccc}
I_{n-1}(3412) &\longrightarrow& I_n(3412) \\
\pi & \mapsto & 1 \oplus \pi
\end{array}
$$
is a bijection between $I_{n-1}(3412)$ and the set of involutions in $I_n(3412)$ which begin with 1.
\item
For all $n \ge 0$ and all $j$ such that $2 \le j \le n$, the map
$$
\begin{array}{ccc}
I_{j-2}(3412) \times I_{n-j}(3412) &\longrightarrow& I_n(3412) \\
(\pi,\sigma) &\mapsto& \pi * \sigma \\
\end{array}
$$
is a bijection between $I_{j-2}(3412) \times I_{n-j}(3412)$ and the set of involutions in $I_n(3412)$ which begin with $j$.
\end{enumerate}
\end{proposition}

\section{Continued Fractions}

In this section we will encounter several continued fractions, for which we will use the following notation.

\begin{definition}
For any given expressions $a_i$ $(i \ge 0)$ and $b_i$ $(i \ge 0)$ we write
$$\frac{a_0}{b_0} \cplus \frac{a_1}{b_1} \cplus \frac{a_2}{b_2} \cplus \frac{a_3}{b_3} \cplus \ldots$$
to denote the infinite continued fraction
$$
\frac{a_0}{\displaystyle
b_0 + \frac{a_1}{\displaystyle
b_1 + \frac{a_2}{\displaystyle
b_2 + \frac{a_3}{\displaystyle
b_3 + \frac{a_4}{\displaystyle
b_4 + \cdots}}}}}.
$$
\end{definition}

For all $i \ge 1$, let $x_i$ denote an indeterminate, let ${\bf x} = \langle x_1, x_2, \ldots \rangle$, and let $C^+({\bf x})$ and $C^-({\bf x})$ denote the generating functions given by
$$C^+({\bf x}) = \sum_{\pi \in I(3412)} \prod_{k \ge 1} x^{\tau_k(\pi)}$$
and
$$C^-({\bf x}) = \sum_{\pi \in I(3412)} (-1)^{sign(\pi)} \prod_{k \ge 1} x^{\tau_k(\pi)}.$$
Egge has shown \cite[Theorem 3.3]{Egge} that
\begin{equation}
\label{eqn:cf}
C^+({\bf x}) =
\frac{1}{1-x_1} \cminus \frac{x_1^2 x_2}{1 - x_1 x_2^2 x_3} \cminus \frac{x_1^2 x_2^5 x_3^4 x_4}{1 - x_1 x_2^4 x_3^6 x_4^4 x_5} \cminus \cdots \cminus \frac{\prod\limits_{i=1}^{2n} x_i^{{{2n-2} \choose {i-1}} + {{2n-1} \choose {i-1}}}}{1 - \prod\limits_{i=1}^{2n+1} x_i^{{{2n} \choose {i-1}}}} \cminus \cdots.
\end{equation}
In this section we refine this result by expressing the signed generating function $C^-({\bf x})$ as a continued fraction.
Combining this expression with (\ref{eqn:cf}) leads to expressions for the restriction of $C^+({\bf x})$ to even (resp. odd) permutations as a sum (resp. difference) of two continued fractions.
We begin with the following simple observation.

\begin{lemma}
For all $i \ge 1$, let $x_i$ denote an indeterminate.
Then
\begin{equation}
\label{eqn:negativex2}
C^-({\bf x}) = C^+(x_1, -x_2, x_3, x_4, \ldots).
\end{equation}
\end{lemma}
\begin{proof}
Fix $\pi \in I(3412)$.
If $\pi$ is even then it contributes $\prod_{k \ge 1} x_k^{\tau_k(\pi)}$ on both sides of (\ref{eqn:negativex2}) and if $\pi$ is odd then it contributes $-\prod_{k \ge 1} x_k^{\tau_k(\pi)}$ on both sides of (\ref{eqn:negativex2}).
\end{proof}

We now express $C^-({\bf x})$ as a continued fraction.

\begin{proposition}
For all $i \ge 1$, let $x_i$ denote an indeterminate.
Then
\begin{equation}
\label{eqn:cfevenminusodd}
C^-({\bf x}) = \frac{1}{1-x_1} \cplus \frac{x_1^2 x_2}{1 - x_1 x_2^2 x_3} \cplus \frac{x_1^2 x_2^5 x_3^4 x_4}{1 - x_1 x_2^4 x_3^6 x_4^4 x_5} \cplus \cdots \cplus \frac{\prod\limits_{i=1}^{2n} x_i^{{{2n-2} \choose {i-1}} + {{2n-1} \choose {i-1}}}}{1 - \prod\limits_{i=1}^{2n+1} x_i^{{{2n} \choose {i-1}}}} \cplus \cdots.
\end{equation}
\end{proposition}
\begin{proof}
Replace $x_2$ with $-x_2$ in (\ref{eqn:cf}) and use (\ref{eqn:negativex2}).
\end{proof}

Using (\ref{eqn:cfevenminusodd}) we can express this signed generating function with respect to various statistics as a continued fraction.

\begin{corollary}
For any permutation $\pi$, let $inv(\pi)$ denote the number of inversions in $\pi$.
Then
$$\sum_{\pi \in I(3412)} (-1)^{sign(\pi)} q^{inv(\pi)} x^{|\pi|} =
\frac{1}{1 - x_1} \cplus \frac{x^2 q}{1 - x q^2} \cplus \frac{x^2 q^5}{1-x q^4} \cplus \cdots \cplus \frac{x^2 q^{4n+1}}{1-x q^{2n+2}} \cplus \cdots.$$
\end{corollary}
\begin{proof}
In (\ref{eqn:cfevenminusodd}), set $x_1 = x, x_2 = q$, and $x_i = 1$ for all $i \ge 3$.
\end{proof}

For our next application of (\ref{eqn:cfevenminusodd}), recall that a {\em left-to-right maximum} in a permutation $\pi$ is an entry of $\pi$ which is greater than all of the entries to its left.
Similarly, a {\em right-to-left minimum} in $\pi$ is an entry of $\pi$ which is less than all of the entries to its right.
We write $lrmax(\pi)$ to denote the number of left-to-right maxima in $\pi$ and we write $rlmin(\pi)$ to denote the number of right-to-left minima in $\pi$.
Egge has shown \cite[Proposition 3.7]{Egge} that if $\pi \in I(3412)$ then
\begin{equation}
\label{eqn:lrmaxrlmin}
lrmax(\pi) = rlmin(\pi) = \sum_{k=1}^\infty (-1)^{k-1} \tau_k(\pi).
\end{equation}
Combining this with (\ref{eqn:cfevenminusodd}) gives us the following continued fraction expansion.

\begin{corollary}
We have
\begin{eqnarray*}
\lefteqn{\sum_{\pi \in I(3412)} (-1)^{sign(\pi)} q^{lrmax(\pi)} x^{|\pi|} = \sum_{\pi \in I(3412)} (-1)^{sign(\pi)} q^{rlmin(\pi)} x^{|\pi|} =} & & \\
 & & \frac{1}{1-x} \cplus \frac{x^2 q}{1-x} \cplus \frac{x^2}{1-x} \cplus \cdots \cplus \frac{x^2}{1-x} \cplus \cdots.
\end{eqnarray*}
and
\begin{displaymath}
\sum_{\pi \in I(3412)} (-1)^{sign(\pi)} q^{lrmax(\pi)} x^{|\pi|} = \sum_{\pi \in I(3412)} (-1)^{sign(\pi)} q^{rlmin(\pi)} x^{|\pi|} = \frac{2}{(2-q)(1-x) - q\sqrt{1-2x+5x^2}}.
\end{displaymath}
\end{corollary}
\begin{proof}
To obtain the first line, set $x_1 = xq$ and $x_i = q^{(-1)^{i-1}}$ for all $i \ge 2$ in (\ref{eqn:cfevenminusodd}) and use (\ref{eqn:lrmaxrlmin}) to simplify the result.
The second line follows routinely from the first line.
\end{proof}

For our final application of (\ref{eqn:cfevenminusodd}), recall that $i$ is a fixed point for a permutation $\pi$ whenever $\pi(i) = i$.
We write $fix(\pi)$ to denote the number of fixed points in $\pi$.
Egge has shown \cite[Proposition 3.9]{Egge} that if $\pi \in I(3412)$ then
\begin{equation}
\label{eqn:fixed}
fix(\pi) = \sum_{k=1}^\infty (-2)^{k-1} \tau_k(\pi).
\end{equation}
Combining this with (\ref{eqn:cfevenminusodd}) gives us the following continued fraction expansion.

\begin{corollary}
We have
$$\sum_{\pi \in I(3412)} (-1)^{sign(\pi)} q^{fix(\pi)} x^{|\pi|} = \frac{1}{1-xq} \cplus \frac{x^2}{1-xq} \cplus \frac{x^2}{1-xq} \cplus \cdots \cplus \frac{x^2}{1-xq} \cplus \cdots$$
and
$$\sum_{\pi \in I(3412)} (-1)^{sign(\pi)} q^{fix(\pi)} x^{|\pi|} = \frac{-1+xq-\sqrt{1-2xq+x^2q^2+4x^2}}{2x^2}.$$
\end{corollary}
\begin{proof}
To obtain the first line, set $x_1 = xq$ and $x_i = q^{(-2)^{i-1}}$ for all $i \ge 2$ in (\ref{eqn:cfevenminusodd}) and use (\ref{eqn:fixed}) to simplify the result.
The second line follows routinely from the first line.
\end{proof}

\section{Involutions Which Avoid 3412 and Another Pattern}

We now turn our attention to signed generating functions for involutions in $I(3412)$ which avoid a set of additional patterns.
We begin by recalling a method of decomposing permutations and a map on permutations.

\begin{definition}
Fix $n \ge 1$.
We call a permutation $\pi \in S_n$ {\em direct sum indecomposable} whenever there do not exist nonempty permutations $\pi_1$ and $\pi_2$ such that $\pi = \pi_1 \oplus \pi_2$.
\end{definition}

Observe that for every permutation $\pi$ there exists a unique sequence of direct sum indecomposable permutations $\alpha_1, \ldots, \alpha_k$ such that $\pi = \alpha_1 \oplus \alpha_2 \oplus \cdots \oplus \alpha_k.$

\begin{definition}
\label{defn:overline}
For any permutation $\pi$, we define $\overline{\pi}$ as follows.
\begin{enumerate}
\item
$\overline{\emptyset} = \emptyset$ and $\overline{1} = \emptyset$.
\item
If $|\pi| \ge 2$ and there exists a permutation $\sigma$ such that $\pi = 1 \ominus \sigma \ominus 1$ then $\overline{\pi} = \sigma$.
\item
If $|\pi| \ge 2$, there exists a permutation $\sigma$ such that $\pi = 1 \ominus \sigma$, and $\sigma$ does not end with 1 then $\overline{\pi} = \sigma$.
\item
If $|\pi| \ge 2$, there exists a permutation $\sigma$ such that $\pi = \sigma \ominus 1$, and $\pi$ does not begin with $|\pi|$ then $\overline{\pi} = \sigma$.
\item
If $|\pi| \ge 2$, $\pi$ does not begin with $|\pi|$, and $\pi$ does not end with 1 then $\overline{\pi} = \pi$.
\end{enumerate}
\end{definition}

\noindent
Observe that if $\pi$ and $\sigma$ are permutations then $1 \ominus \pi \ominus 1$ avoids $\sigma$ if and only if $\pi$ avoids $\overline{\sigma}$.

For any set $T$ of permutations we write
$$F^+_T(x) = \sum_{\pi \in I(3412, T)} x^{|\pi|},$$
and
$$F^-_T(x) = \sum_{\pi \in I(3412, T)} (-1)^{sign(\pi)} x^{|\pi|}.$$
Egge has shown \cite[Corollary 5.6]{Egge} that if $\pi = \alpha_1 \oplus \cdots \oplus \alpha_k$ is a permutation and $\alpha_1,\ldots,\alpha_k$ are direct sum indecomposable then
\begin{equation}
\label{eqn:Fpiplusrecurrence}
F^+_\pi(x) = 1 + x F^+_\beta(x) + x^2 \sum_{i=1}^k \left( F^+_{\overline{\alpha_1 \oplus \cdots \oplus \alpha_i}}(x) - F^+_{\overline{\alpha_1 \oplus \cdots \oplus \alpha_{i-1}}}(x)\right) F^+_{\alpha_i \oplus \cdots \oplus \alpha_k}(x),
\end{equation}
where $\beta = \pi$ if $\alpha_1 \neq 1$ and $\beta = \alpha_2 \oplus \cdots \oplus \alpha_k$ if $\alpha_1 = 1$.
Our main result in this section is a similar recurrence relation for $F^-_\pi(x)$.
Combining this with (\ref{eqn:Fpiplusrecurrence}) allows one to compute the generating function for the even (or odd) involutions in $I(3412,\pi)$ for any permutation $\pi$.

\begin{theorem}
\label{thm:Fpiminusrecurrence}
Suppose $\pi = \alpha_1 \oplus \cdots \oplus \alpha_k$ is a permutation, where $\alpha_1,\ldots,\alpha_k$ are direct sum indecomposable.
Then
\begin{equation}
\label{eqn:Fpiminusrecurrence}
F^-_\pi(x) = 1 + x F^-_\beta(x) - x^2 \sum_{i=1}^k \left( F^-_{\overline{\alpha_1 \oplus \cdots  \oplus \alpha_i}}(x) - F^-_{\overline{\alpha_1 \oplus \ldots \oplus \alpha_{i-1}}}(x)\right) F^-_{\alpha_i \oplus \cdots \oplus \alpha_k}(x).
\end{equation}
Here $\beta = \pi$ if $\alpha_1 \neq 1$ and $\beta = \alpha_2 \oplus \cdots \oplus \alpha_k$ if $\alpha_1 = 1$.
\end{theorem}
\begin{proof}
The set $I(3412,\pi)$ can be partitioned into three sets:  the set $A_1$ containing only the empty permutation, the set $A_2$ of those involutions which begin with 1, and the set $A_3$ of those involutions which do not begin with 1.

The set $A_1$ contributes 1 to the desired generating function.

In view of Proposition \ref{prop:permmap}(i), the set $A_2$ contributes $x F^-_\beta(x)$ to the desired generating function, where $\beta = \pi$ if $\alpha_1 \neq 1$ and $\beta = \alpha_2 \oplus \cdots \oplus \alpha_k$ if $\alpha_1 = 1$.

To obtain the contribution of $A_3$ to the desired generating function, we first observe that in view of Proposition \ref{prop:permmap}(ii), all permutations in $A_3$ have the form $\sigma_1 * \sigma_2$.
Since each $\alpha_i$ is direct sum indecomposable, if $\sigma_1 * \sigma_2$ contains a subsequence of type $\alpha_i$ then that subsequence is entirely contained in either $1 \ominus \sigma_1 \ominus 1$ or $\sigma_2$.
As a result, the set of involutions which avoid 3412 and $\pi$ and which do not begin with 1 can be partitioned into sets $B_1, \ldots, B_k$, where $B_i$ is the set of such involutions in which $\sigma_1$ contains $\overline{\alpha_1 \oplus \cdots \oplus \alpha_{i-1}}$ but avoids $\overline{\alpha_1 \oplus \cdots \oplus \alpha_i}$.
Since $I_n(3412, \overline{\alpha_1 \oplus \cdots \oplus \alpha_{i-1}}) \subseteq I_n(3412, \overline{\alpha_1 \oplus \cdots \oplus \alpha_i})$, and since the sign of $\sigma_1 * \sigma_2$ is the negative of the product of the signs of $\sigma_1$ and $\sigma_2$, the contribution of the set $A_3$ to the desired generating function is
$$- x^2 \sum_{i=1}^k \left( F^-_{\overline{\alpha_1 \oplus \cdots \oplus \alpha_i}}(x) - F^-_{\overline{\alpha_1 \oplus \ldots \oplus \alpha_{i-1}}}(x)\right) F^-_{\alpha_i \oplus \cdots \oplus \alpha_k}(x).$$
Add the contributions of $A_1$, $A_2$, and $A_3$ to obtain (\ref{eqn:Fpiminusrecurrence}).
\end{proof}

Theorem \ref{thm:Fpiminusrecurrence} above is an analogue of \cite[Corollary 5.6]{Egge}, which is a special case of \cite[Theorem 5.5]{Egge}.
Using Egge's techniques, one can also prove an analogue of \cite[Theorem 5.5]{Egge}.
To state this result, we first set some notation.

\begin{definition}
\label{defn:setnotation}
Let $T = \{\pi_1,\ldots,\pi_m\}$ denote a set of permutations and fix direct sum indecomposable permutations $\alpha^i_j$, $1 \le i \le m$, $1 \le j \le k_i$, such that $\pi_i = \alpha^i_1 \oplus \cdots \oplus \alpha^i_{k_i}$.
For all $i_1,\ldots, i_m$ such that $0 \le i_j \le k_j$, let $T^{right}_{i_1,\ldots,i_m} = \{\alpha^1_{i_1} \oplus \cdots \oplus \alpha^1_{k_1},\ldots,\alpha^m_{i_m} \oplus \cdots \oplus \alpha^m_{k_m}\}$.
For any subset $Y \subseteq \{1,\ldots,m\}$, set
\begin{displaymath}
T_Y = \bigcup_{j \in Y} \{\overline{\alpha^j_1 \oplus \cdots \oplus \alpha^j_{i_j-1}}\} \bigcup_{j \not\in Y, 1 \le j \le m} \{\overline{\alpha^j_1 \oplus \cdots \oplus \alpha^j_{i_j}}\}.
\end{displaymath}
\end{definition}

\begin{theorem}
\label{thm:FTminusrecurrence}
With reference to Definition \ref{defn:setnotation},
\begin{displaymath}
F_T^-(x) = 1 + x F^-_{\beta(T)}(x) - x^2 \sum_{i_1,\ldots,i_m = 1}^{k_1,\ldots,k_m} \left( \sum_{Y \subseteq \{1,2,\ldots,m\}} (-1)^{|Y|} F^-_{T_Y}(x) \right) F^-_{T^{right}_{i_1,\ldots,i_m}}(x).
\end{displaymath}
Here $\beta(\pi_i) = \pi_i$ if $\alpha^i_1 \neq 1$, $\beta(\pi_i) = \alpha^i_2 \oplus \cdots \oplus \alpha_{k_i}^i$ if $\alpha^i_1 = 1$, and $\beta(T)$ is the set of permutations obtained by applying $\beta$ to every element of $T$.
\end{theorem}

We omit the proof of Theorem \ref{thm:FTminusrecurrence} for the sake of brevity.

For the remainder of this section we use (\ref{eqn:Fpiminusrecurrence}) to find $F^-_\pi(x)$ for various $\pi$.
Combining these results with Egge's expressions for various $F_\pi^+(x)$ allows one to find the generating function for the even (or odd) involutions in $I(3412,\pi)$ for these $\pi$.
We express all of our generating functions in terms of Chebyshev polynomials of the second kind, so we begin by recalling these polynomials.

\begin{definition}
For all $n$ we write $U_n(x)$ to denote the {\em $n$th Chebyshev polynomial of the second kind}, which is defined by $U_n(x) = 0$ for $n < 0$ and ${\displaystyle U_n(\cos t) = \frac{\sin ((n+1)t)}{\sin t}}$ for $n \ge 0$.
These polynomials satisfy
\begin{equation}
\label{eqn:Chebyshevrecurrence}
U_n(x) = 2x U_{n-1}(x) - U_{n-2}(x) \hspace{30pt} (n \neq 0).
\end{equation}
\end{definition}

We will often use two specializations of $U_n(x)$, which are defined for all $n$ by
$$V_n = V_n(x) = U_n\left(\frac{1-x}{2x}\right)$$
and
$$W_n = W_n(x) = U_n\left(\frac{1-x}{2ix}\right).$$
Observe that by (\ref{eqn:Chebyshevrecurrence}) we have
\begin{displaymath}
x V_n(x) = (1-x) V_{n-1}(x) - x V_{n-2}(x)
\end{displaymath}
and
\begin{equation}
\label{eqn:ourChebyWrecurrence}
ix W_n(x) = (1-x) W_{n-1}(x) - ix W_{n-2}(x).
\end{equation}

Turning our attention to $F^+_\pi(x)$ and $F^-_\pi(x)$, we now consider permutations of the forms $[j] \ominus \pi \ominus [j]$, $[j] \ominus \pi$, and $\pi \ominus [j]$.

\begin{proposition}
\label{prop:Fkpi1minus}
For any permutation $\pi$ we have
\begin{equation}
\label{eqn:Fkpi1plus}
F^+_{1 \ominus \pi \ominus 1}(x) = \frac{1}{1-x-x^2 F^+_\pi(x)}.
\end{equation}
and
\begin{equation}
\label{eqn:Fkpi1minus}
F_{1 \ominus \pi \ominus 1}^-(x) = \frac{1}{1-x+x^2 F_\pi^-(x)}.
\end{equation}
Moreover, if $F_\pi^+(x) = \frac{f_0(x)}{f_1(x)}$ for polynomials $f_0$ and $f_1$ then for all $j \ge 1$,
\begin{equation}
\label{eqn:Fkjpi1jplus}
F^+_{[j] \ominus \pi \ominus [j]}(x) = \frac{f_1(x) V_{j-1}(x) - x f_0(x) V_{j-2}(x)}{x f_1(x) V_j(x) - x^2 f_0(x) V_{j-1}(x)}.
\end{equation}
Similarly, if $F_\pi^-(x) = \frac{f_0(x)}{f_1(x)}$ for polynomials $f_0$ and $f_1$ then for all $j \ge 1$,
\begin{equation}
\label{eqn:Fkjpi1jminus}
F_{[j] \ominus \pi \ominus [j]}^-(x) = \frac{f_1(x) W_{j-1}(x) - ix f_0(x) W_{j-2}(x)}{ix f_1(x) W_j(x) + x^2 f_0(x) W_{j-1}(x)}.
\end{equation}
\end{proposition}
\begin{proof}
To prove (\ref{eqn:Fkpi1minus}), replace $\pi$ with $1 \ominus \pi \ominus 1$ in (\ref{eqn:Fpiminusrecurrence}), use the fact that $1 \ominus \pi \ominus 1$ is direct sum indecomposable and $\beta = 1 \ominus \pi \ominus 1$ to simplify the result, and solve for $F^-_{1 \ominus \pi \ominus 1}(x)$.

To prove (\ref{eqn:Fkjpi1jminus}) we argue by induction on $j$.
When $j = 1$ the result is immediate from (\ref{eqn:Fkpi1minus}), so we assume $j \ge 2$ and the result holds for $j-1$.
Replace $\pi$ with $[j-1] \ominus \pi \ominus [j-1]$ in (\ref{eqn:Fkpi1minus}), use induction to eliminate $F^-_{[j-1] \ominus \pi \ominus [j-1]}(x)$ on the right, and use (\ref{eqn:ourChebyWrecurrence}) to simplify the result and obtain (\ref{eqn:Fkjpi1jminus}).

The proofs of (\ref{eqn:Fkpi1plus}) and (\ref{eqn:Fkjpi1jplus}) are similar to the proofs of (\ref{eqn:Fkpi1minus}) and (\ref{eqn:Fkjpi1jminus}) respectively.
\end{proof}

Arguing as in the proof of Proposition \ref{prop:Fkpi1minus}, one can show that if $\pi$ does not end with 1 then (\ref{eqn:Fkpi1plus})--(\ref{eqn:Fkjpi1jminus}) hold when $1 \ominus \pi \ominus 1$ (resp. $[j]\ominus \pi \ominus [j]$) is replaced with $1 \ominus \pi$ (resp. $[j] \ominus \pi$).
Similarly, one can show that if $\pi$ does not begin with $|\pi|$ then (\ref{eqn:Fkpi1plus})--(\ref{eqn:Fkjpi1jminus}) hold when $1 \ominus \pi \ominus 1$ (resp. $[j]\ominus \pi \ominus [j]$) is replaced with $\pi \ominus 1$ (resp. $\pi \ominus [j]$).

With Proposition \ref{prop:Fkpi1minus} and its analogues in hand, we are ready to compute $F^-_\pi(x)$ for various $\pi$.
We begin with $\pi = [k]$.

\begin{theorem}
\label{thm:F321}
For all $k \ge 1$ we have
\begin{equation}
\label{eqn:F2kminus}
F^-_{[2k]}(x) = \frac{W_{k-1}(x)}{ix W_k(x)}
\end{equation}
and
\begin{equation}
\label{eqn:F2k-1minus}
F^-_{[2k-1]}(x) = \frac{W_{k-2}(x) - ix W_{k-3}(x)}{ix\left( W_{k-1}(x) - ix W_{k-2}(x)\right)}.
\end{equation}
\end{theorem}
\begin{proof}
To prove (\ref{eqn:F2kminus}), set $\pi = \emptyset$ and $j = k$ in (\ref{eqn:Fkjpi1jplus}) and use the fact that $F^-_\pi(x) = \frac{0}{1}$.

The proof of (\ref{eqn:F2k-1minus}) is similar to the proof of (\ref{eqn:F2kminus}).
\end{proof}

\noindent
{\bf Remark}
Lines (\ref{eqn:F2kminus}) and (\ref{eqn:F2k-1minus}) can also be obtained using the methods of Section \ref{sec:rk21} below.

\medskip

Using the analogues of Proposition \ref{prop:Fkpi1minus} along with the fact that $F_{12}^-(x) = \frac{1+x}{1+x^2}$, we obtain the following result.

\begin{theorem}
\label{thm:F231312}
For all $k \ge 1$ we have
\begin{equation}
\label{eqn:Fk312}
F^-_{[k] \ominus 231}(x) = F^-_{[k] \ominus 12}(x) = \frac{(1+x^2) W_{k-1}(x) - ix(1+x) W_{k-2}(x)}{ix \left((1+x^2) W_k(x)-ix(1+x)W_{k-1}(x)\right)}.
\end{equation}
\end{theorem}

\noindent
{\bf Remark}
The fact that $F^-_{[k] \ominus 231}(x) = F^-_{[k] \ominus 12}(x)$ is immediate from \cite[Theorem 6.3]{Egge}, which implies that $I_n(3412, [k] \ominus 231) = I_n(3412, [k] \ominus 12)$ for all $n \ge 0$ and all $k \ge 3$.

\medskip

We now turn our attention to $F^-_{[k] \ominus 213}(x)$, $F^-_{[k] \ominus 132}(x)$, and $F^-_{[k] \ominus 123}(x)$.

\begin{theorem}
\label{thm:F213132}
For all $k \ge 0$ we have
\begin{equation}
\label{eqn:F213132}
F^-_{[k] \ominus 213}(x) = F^-_{[k] \ominus 132}(x) = \frac{W_k(x) - ix W_{k-1}(x)}{ix \left( W_{k+1}(x) - ix W_k(x)\right)}.
\end{equation}
\end{theorem}
\begin{proof}
First observe that if we apply the reverse complement map to $[k] \ominus 213$ and take the inverse of the result we obtain $[k] \ominus 132$, and that these operations preserve parity, so $F^-_{[k] \ominus 213}(x) = F^-_{[k] \ominus 132}(x)$.

To show that $F^-_{[k] \ominus 132}(x)$ is equal to the quantity on the right, first set $\pi = 132$ in (\ref{eqn:Fpiminusrecurrence}), use the fact that $F^-_{21}(x) = \frac{1}{1-x}$, and solve the resulting equation to find that $F_{132}^-(x) = \frac{1}{1-x+x^2}$.
Therefore (\ref{eqn:F213132}) holds for $k = 0$.
To see that the result holds for $k > 0$, first set $j = k$ and $\pi = 132$ in the analogue of (\ref{eqn:Fkjpi1jminus}) in which $[j] \ominus \pi \ominus [j]$ is replaced with $[j] \ominus \pi$.
Then use (\ref{eqn:ourChebyWrecurrence}) and the fact that $F^-_{132}(x) = \frac{1}{1-x+x^2}$ to simplify the result.
\end{proof}

\begin{theorem}
\label{thm:F123}
For all $k \ge 0$ we have
$$F^-_{[k] \ominus 123}(x) = \frac{(1+x^2)^3 W_{k-1}(x) - ix (1+3x^2+2x^4+x^5) W_{k-2}(x)}{ix\left((1+x^2)^3 W_k(x) - ix (1+3x^2+2x^4+x^5) W_{k-1}(x)\right)}.$$
\end{theorem}
\begin{proof}
This is similar to the second half of the proof of Theorem \ref{thm:F213132}.
\end{proof}

Recall that 213 and 123 are examples of layered permutations, which are defined as follows.

\begin{definition}
Fix $n \ge 1$ and let $l_1, l_2, \ldots, l_m$ denote a sequence such that $l_i \ge 1$ for $1 \le i \le m$ and $\sum\limits_{i=1}^m l_i = n$.
We write $[l_1, l_2, \ldots, l_m]$ to denote the permutation in $S_n$ given by
$$[l_1, l_2, \ldots, l_m] = [l_1] \oplus \cdots \oplus [l_m].$$
We call a permutation {\em layered} whenever it has the form $[l_1,\ldots,l_m]$ for some sequence $l_1, \ldots, l_m$.
\end{definition}

Observe that if $m \ge 2$ then $\overline{[l_1, \ldots, l_m]} = [l_1,\ldots,l_m]$.
In view of (\ref{eqn:Fpiminusrecurrence}), (\ref{eqn:F2kminus}), and (\ref{eqn:F2k-1minus}), the generating function $F^-_{[l_1, \ldots, l_m]}(x)$ can be expressed in terms of Chebyshev polynomials of the second kind for any layered permutation $[l_1,\ldots, l_m]$.
To do this for $m=2$ we will use the following well-known identity for Chebyshev polynomials.

\begin{lemma}
\label{lem:Chebysum}
For all $k, l \ge -1$ and all $w \ge 0$ we have
\begin{equation}
\label{eqn:Chebysum}
U_{k+w} U_{l+w} - U_k U_l = U_{w-1} U_{k+l+w+1}.
\end{equation}
\end{lemma}

It will also be useful to record the recurrence relations for $F^-_{[k,l]}(x)$.

\begin{lemma}
For all $k \ge 2$ and all $l \ge 1$ we have
\begin{equation}
\label{eqn:F[1l]}
F^-_{[1,l]}(x) = \frac{1 + x F^-_{[l]}(x)}{1+x^2 F^-_{[l]}(x)}
\end{equation}
and
\begin{equation}
\label{eqn:F[kl]}
F^-_{[k,l]}(x) = \frac{1+x^2 F^-_{[k-2]}(x) F^-_{[l]}(x)}{1-x+x^2 F^-_{[k-2]}(x) + x^2 F^-_{[l]}(x)}.
\end{equation}
\end{lemma}
\begin{proof}
To prove (\ref{eqn:F[1l]}), set $\pi = [1,l]$ in (\ref{eqn:Fpiminusrecurrence}) and solve the resulting equation for $F^-_{[1,l]}(x)$.

The proof of (\ref{eqn:F[kl]}) is similar to the proof of (\ref{eqn:F[1l]}).
\end{proof}

We now compute $F^-_{[k,l]}(x)$ when $k$ and $l$ are not both odd.

\begin{theorem}
\label{thm:2layered}
For all $k, l \ge 1$ such that $k$ and $l$ are not both odd we have
\begin{equation}
\label{eqn:2layeredkl}
F^-_{[k,l]}(x) = F^-_{[k+l]}(x).
\end{equation}
\end{theorem}
\begin{proof}
We consider four cases:  $k = 1$ and $l$ is even, $k$ and $l$ are both even, $k$ is even and $l$ is odd, and $k > 1$ is odd and $l$ is even.
All four cases are similar, so we only give the details for the case in which $k >1$ is even and $l$ is odd.

In (\ref{eqn:F[kl]}) replace $k$ with $2k$ and $l$ with $2l-1$, use (\ref{eqn:F2kminus}) and (\ref{eqn:F2k-1minus}) to write the result in terms of Chebyshev polynomials, clear denominators, and use (\ref{eqn:ourChebyWrecurrence}) in the resulting denominator to obtain
\begin{displaymath}
F^-_{[2k,2l-1]}(x) = \frac{W_{k-1} \left( W_{l-1} - ix W_{l-2} \right) - W_{k-2}\left( W_{l-2} - ix W_{l-3}\right)}{ix\left( W_k \left( W_{l-1} - ix W_{l-2} \right) - W_{k-1} \left( W_{l-2} - ix W_{l-3} \right)\right)}.
\end{displaymath}
Now apply (\ref{eqn:Chebysum}) to the numerator and denominator to obtain the right side of (\ref{eqn:F2k-1minus}) with $k$ replaced by $k+l$, as desired.
\end{proof}

Next we compute $F_{[k,l]}(x)$ when $k$ and $l$ are both odd.

\begin{theorem}
For all $k,l \ge 1$ we have
\begin{equation}
\label{eqn:2layeredoddodd}
F^-_{[2k-1,2l-1]}(x) = \frac{W_{k+l-1}(x) - 2i W_{k+l-2}(x)-W_{k+l-3}(x)}{ix( W_{k+l}(x) - 2i W_{k+l-1}(x) - W_{k+l-2}(x))}.
\end{equation}
\end{theorem}
\begin{proof}
We consider two cases:  $k = 1$ and $k > 1$.
These are similar, so we only give details for the case in which $k > 1$.

In (\ref{eqn:F[kl]}) replace $k$ with $2k-1$ and $l$ with $2l-1$, use (\ref{eqn:F2k-1minus}) and (\ref{eqn:ourChebyWrecurrence}) to write the result in terms of Chebyshev polynomials, and clear denominators to obtain
\begin{displaymath}
F^-_{[2k-1,2l-1]}(x) = \frac{(W_{k-1} - i W_{k-2})(W_l-iW_{l-1}) - (W_{k-2}-iW_{k-3})(W_{l-1} - iW_{l-2})}{(1-x) Z_1-ix Z_2-ix Z_3}.
\end{displaymath}
Here
$$Z_1 = (W_{k-1}-iW_{k-2})(W_l-iW_{l-1}),$$
$$Z_2 = (W_{k-2}-iW_{k-3})(W_l-iW_{l-1}),$$
and
$$Z_3 = (W_{l-1}-iW_{l-2})(W_{k-1}-iW_{k-2}).$$
Group terms in the numerator and denominator according to the power of $i$ contributed by factors of the form $i W_n$ to obtain
$$F^-_{[2k-1,2l-1]}(x) = \frac{B_1 - i B_2 + i^2 B_3}{C_1 - i C_2 + i^2 C_3},$$
where
$$B_1 = W_{k-1} W_l-W_{k-2} W_{l-1},$$
$$B_2 = W_{k-2} W_l + W_{k-1} W_{l-1} - W_{k-2} W_{l-2} - W_{k-3}W_{l-1},$$
$$B_3 = W_{k-2}W_{l-1}-W_{k-3}W_{l-2},$$
$$C_1 = (1-x)W_{k-1}W_l - ixW_{k-2} W_l - ixW_{l-1}W_{k-1},$$
$$C_2 = (1-x)(W_{k-2} W_l + W_{k-1} W_{l-1}) - ix (W_{k-2} W_{l-1} + W_{k-3}W_l) -ix(W_{l-1}W_{k-2} + W_{l-2}W_{k-1}),$$
and
$$C_3 = (1-x) W_{k-2} W_{l-1} - ix W_{k-3} W_{l-1} - ix W_{l-2} W_{k-2}.$$
Apply (\ref{eqn:Chebysum}) to $B_1$ to find $B_1 = W_{k+l-1}$.
Similarly, $B_2 = 2i W_{k+l-2}$ and $B_3 = W_{k+l-3}$.
Now apply (\ref{eqn:ourChebyWrecurrence}) and (\ref{eqn:Chebysum}) to $C_1$ to find $C_1 = ix W_{k+l}$.
Similarly, $C_2 = ix 2i W_{k+l-1}$ and $C_3 = ix W_{k+l-2}$.
Combine these results to obtain (\ref{eqn:2layeredoddodd}), as desired.
\end{proof}

When $m \ge 3$ the generating function $F^-_{[l_1,\ldots,l_m]}(x)$ does not reduce quite as nicely as it does when $m=2$.
For example, using the same techniques as in the proof of Theorem \ref{thm:2layered} one can prove that for all $k_1, k_2, k_3 \ge 1$,
\begin{equation}
\label{eqn:k1k2k3}
F^-_{[2k_1,2k_2,2k_3]}(x) = \frac{W_{k_1+k_2+k_3}(x) W_{k_1+k_2+k_3-1}(x) + W_{k_1+k_2-1}(x) W_{k_1+k_3-1}(x) W_{k_2+k_3-1}(x)}{ix W_{k_1+k_2}(x) W_{k_1+k_3}(x) W_{k_2+k_3}(x)}.
\end{equation}
Nevertheless, (\ref{eqn:2layeredkl}), (\ref{eqn:2layeredoddodd}), and (\ref{eqn:k1k2k3}) suggest the following conjecture.

\begin{conjecture}
\label{conj:minussymmetry}
For all $m \ge 1$ and all $l_1, \ldots, l_m \ge 1$, the generating function $F^-_{[l_1,\ldots, l_m]}(x)$ is symmetric in $l_1, \ldots, l_m$.
\end{conjecture}

We have verified this conjecture in the case $m = 3$ for $l_i \le 24$ and in the case $m = 4$ for $l_i \le 20$ using a Maple program.

Egge has conjectured \cite[Conjecture 6.9]{Egge} that $F^+_{[l_1,\ldots,l_m]}(x)$ is also symmetric in $l_1,\ldots, l_m$.
This suggests the following conjecture.

\begin{conjecture}
\label{conj:evenoddsymmetry}
For all $m \ge 1$ and all $l_1, \ldots, l_m \ge 1$, the generating functions for the even involutions in $I(3412, [l_1, \ldots, l_m])$ and for the odd involutions in $I(3412, [l_1, \ldots, l_m])$ are symmetric in $l_1, \ldots, l_m$.
\end{conjecture}

\section{Involutions Avoiding 3412 and Containing $[k]$}
\label{sec:rk21}

In \cite{Egge} Egge finds the generating function for the involutions in $I(3412)$ which contain exactly $r$ subsequences of type $[k]$.
In this section we prove the following refinements of this result, which give the signed analogue of Egge's generating function.
Throughout we adopt the convention that ${{a} \choose {0}} = 1$ and ${{a} \choose {-1}} = 0$ for any integer $a$.

\begin{theorem}
\label{thm:biggfeven}
Fix $r \ge 1$, $k \ge 1$, and $b \ge 0$ such that
$$r < min\left( {{2k+2b+2} \choose {2k-1}}, {{2k+2b}\choose{2k-1}} + {{2k+2b+1}\choose{2k-1}}\right).$$
Then
\begin{eqnarray*}
\lefteqn{\sum_{\pi \in I(3412)} (-1)^{sign(\pi)} x^{|\pi|} =} & & \\
& & \sum \prod_{s=0}^b {{d_s + d_{s+1} + l_s - 1}\choose{d_{s+1}+l_s}} {{d_{s+1}+l_s}\choose{l_s}} \frac{W_{k-1}^{d_0-1}}{W_k^{d_0+1}} i^{\sum\limits_{j=0}^b 2d_j + d_0-1} x^{-1-d_0+\sum\limits_{j=0}^b (2 d_j + l_j)}.
\end{eqnarray*}
Here the sum on the left is over all involutions in $I(3412)$ which contain exactly $r$ subsequences of type $[2k]$ and the sum on the right is over all sequences $d_0, \ldots, d_b$ and $l_0, \ldots, l_b$ of nonnegative integers such that
\begin{equation}
\label{eqn:bigreven}
r = \sum_{j=0}^b d_j \left( {{2k+2j-2}\choose{2k-1}} + {{2k+2j-1}\choose{2k-1}}\right) + \sum_{j=0}^b l_j {{2k+2j}\choose{2k-1}}.
\end{equation}
\end{theorem}

\begin{theorem}
\label{thm:biggfodd}
Fix $r \ge 1$, $k \ge 1$, and $b \ge 0$ such that ${{2k+2b}\choose{2k}} \le r < {{2k+2b+2}\choose{2k}}$.
Then
\begin{eqnarray*}
\lefteqn{\sum_{\pi \in I(3412)} (-1)^{sign(\pi)} x^{|\pi|} =} & & \\
& & \sum \prod_{s=0}^b {{d_s+d_{s-1}+l_s-1}\choose{d_s+l_s}}{{d_s+l_s}\choose{l_s}} \frac{W_k^{d_0+l_0-1}}{\left(W_{k+1} + W_k\right)^{d_0+l_0+1}} i^{-1-d_0+\sum\limits_{j=0}^b 2d_j} x^{-1-d_0-l_0+\sum\limits_{j=0}^b (2d_j+l_j)}.
\end{eqnarray*}
Here the sum on the left is over all involutions in $I(3412)$ which contain exactly $r$ subsequences of type $[2k]$ and the sum on the right is over all sequences $d_0, \ldots, d_b$ and $l_0, \ldots, l_b$ of nonnegative integers such that
\begin{displaymath}
r = \sum_{j=0}^b d_j\left( {{2k+2j+1}\choose{2k}} + {{2k+2j}\choose{2k}}\right) + \sum_{j=0}^b l_j{{2k+2j}\choose{2k}}.
\end{displaymath}
For notational convenience we set $d_{-1} = 1$.
\end{theorem}

In order to prove these results, we first recall the relationship between $I_n(3412)$ and the set $\cM_n$ of Motzkin paths of length $n$.
These are the lattice paths from $(0,0)$ to $(n,0)$ which contain only up $(1,1)$, down $(1,-1)$, and level $(1,0)$ steps and which do not pass below the line $y = 0$.
We begin by defining a family of statistics on the set $\cM$ of all Motzkin paths.

\begin{definition}
Suppose $\pi$ is a Motzkin path.
For any step $s \in \pi$ we write $ht(s)$ to denote the {\em height} of $s$, which is the $x$-coordinate of the left-most point of $s$.
For any $k$ we write
\begin{equation}
\label{eqn:tauk}
\tau_k(\pi) = \sum_{s \in \pi} {{2 ht(s)} \choose {k-1}} + \sum_{s \in \pi} {{2 ht(s) - 1} \choose {k-1}},
\end{equation}
where the first sum on the right is over all up and level steps in $\pi$ and the second sum on the right is over all down steps in $\pi$.
Here we use the convention that ${{n} \choose {k}} = 0$ whenever $k < 0$ or $k > n$.
\end{definition}

Next we recall a $\tau_k$-preserving bijection between $I_n(3412)$ and $\cM_n$.

\begin{definition}
For any $\pi \in \cM_n$, we write $\varphi(\pi)$ to denote the permutation obtained as follows.
Number the steps in $\pi$ from left to right with $1, 2, \ldots, n$.
For each up step at height $k$, find the first down step at height $k+1$ to its right and switch the labels of the two steps.
Then $\varphi(\pi)$ is the involution obtained by reading the resulting labels from left to right.
\end{definition}

\begin{proposition}
\label{prop:varphibijection}
(\cite[Proposition 2.13]{Egge})
For all $n \ge 0$, the map $\varphi$ is a bijection between $I_n(3412)$ and $\cM_n$ such that $\tau_k(\varphi(\pi)) = \tau_k(\pi)$ for all $k$ and all $\pi \in \cM$.
\end{proposition}

The bijection $\varphi$ allows us to transfer the sign function from $I_n(3412)$ to $\cM_n$.

\begin{proposition}
For any Motzkin path $\pi$ we write $U(\pi)$ to denote the number of up steps in $\pi$ and we write $D(\pi)$ to denote the number of down steps in $\pi$.
Then for any $\pi \in \cM$,
\begin{equation}
\label{eqn:invUD}
i^{U(\pi)+D(\pi)} = sign(\varphi(\pi)).
\end{equation}
\end{proposition}
\begin{proof}
Set $k = 2$ in (\ref{eqn:tauk}) and use Proposition \ref{prop:varphibijection} to find that $(-1)^{D(\pi)} = (-1)^{\tau_2(\varphi(\pi))}$ for all $\pi \in \cM$.
Now the result follows, since $D(\pi) = U(\pi)$ for all $\pi \in \cM$ and $sign(\pi) = (-1)^{\tau_2(\varphi(\pi))}$ for all $\pi \in I(3412)$.
\end{proof}

We now turn our attention to three families of matrices which will use in proving our main results.

\begin{definition}
For all $k \ge 0$ we write $A_k$ to denote the $k+1$ by $k+1$ tridiagonal matrix given by
$$A_k = \left( \matrix{x & ix & 0 & 0 & 0 & \cdots & 0 & 0 & 0 & 0 \cr
ix & x & ix & 0 & 0 & \cdots & 0 & 0 & 0 & 0 \cr
0 & ix & x & ix & 0 & \cdots & 0 & 0 & 0 & 0 \cr
\vdots & & \ddots & \ddots & \ddots & & & & \vdots & \vdots \cr
\vdots & & & \ddots & \ddots & \ddots & & & \vdots & \vdots \cr
\vdots & & & & \ddots & \ddots & \ddots & & \vdots & \vdots \cr
\vdots & & & & & \ddots & \ddots & \ddots & \vdots & \vdots \cr
0 & 0 & 0 & 0 & 0 & \cdots & ix & x & ix & 0 \cr
0 & 0 & 0 & 0 & 0 & \cdots & 0 & ix & x & ix \cr
0 & 0 & 0 & 0 & 0 & \cdots & 0 & 0 & ix & x}
\right).$$
We write $B_k$ to denote the $k+1$ by $k+1$ tridiagonal matrix obtained by replacing the entry in the lower right corner of $A_k$ with 0.
We write $C_k$ to denote the $k+1$ by $k+1$ tridiagonal matrix obtained by replacing the entry in the upper left corner of $A_k$ with 0.
\end{definition}

The matrices $A_k$, $B_k$, and $C_k$ are closely related to generating functions for various sets of Motzkin paths.
To describe this relationship, we let $\cM(r,s,k)$ denote the set of lattice paths involving only up $(1,1)$, down $(1,-1)$, and level $(1,0)$ steps which begin at a point at height $r$, $0 \le r \le k$, end at a point at height $s$, $0 \le s \le k$, and do not cross the lines $y = k$ and $y = 0$.
Similarly, we let $\cN(r,s,k)$ denote the set of lattice paths in $\cM(r,s,k)$ which do not have any level steps at height $k$, and we let $\cO(r,s,k)$ denote the set of lattice paths in $\cM(r,s,k)$ which do not have any level steps at height 0.
Modifying the proof of \cite[Theorem A2]{Krattenthaler} slightly, we find that
\begin{equation}
\label{eqn:Agf}
\sum_{\pi \in \cM(r,s,k)} i^{U(\pi) + D(\pi)} x^{|\pi|} = \frac{(-1)^{r+s} \det(I - A_k; s, r)}{\det(I - A_k)},
\end{equation}
\begin{equation}
\label{eqn:Bgf}
\sum_{\pi \in \cN(r,s,k)} i^{U(\pi) + D(\pi)} x^{|\pi|} = \frac{(-1)^{r+s} \det(I - B_k; s, r)}{\det(I - B_k)},
\end{equation}
and
\begin{equation}
\label{eqn:Cgf}
\sum_{\pi \in \cO(r,s,k)} i^{U(\pi) + D(\pi)} x^{|\pi|} = \frac{(-1)^{r+s} \det(I - C_k; s, r)}{\det(I - C_k)}.
\end{equation}
Here $|\pi|$ is the number of steps in $\pi$, the quantity $U(\pi)$ (resp. $D(\pi))$ is the number of up (resp. down) steps in $\pi$, the matrix $I$ is the identity matrix of the appropriate size, and $\det(M; s, r)$ is the minor of the matrix $M$ in which the $s$th row and $r$th column of $M$ have been deleted.
We will use the fact that the determinants in (\ref{eqn:Agf}), (\ref{eqn:Bgf}), and (\ref{eqn:Cgf}) can be expressed in terms of Chebyshev polynomials of the second kind.
In particular, arguing by induction on $k$ we find that for all $k \ge 0$,
\begin{equation}
\label{eqn:Auk}
(ix)^{k+1} W_{k+1}(x) = \det(I - A_k)
\end{equation}
and
\begin{equation}
\label{eqn:Buk}
(ix)^{k+1} \left( W_{k+1}(x) + W_k(x) \right) = \det(I - B_k) = \det(I-C_k).
\end{equation}

We now prove Theorem \ref{thm:biggfeven}.

\noindent
{\em Proof of Theorem \ref{thm:biggfeven}}.
For any Motzkin path $\pi$, let the weight of each up or down step in $\pi$ be $ix$, let the weight of each level step be $x$, and let the weight of $\pi$ be the product of the weights of its steps.
Observe that in view of (\ref{eqn:invUD}) and Proposition \ref{prop:varphibijection}, the desired generating function is the sum of the weights of the Motzkin paths from $(0,0)$ to $(n,0)$ for which $\tau_{2k}(\pi) = r$.
To compute this generating function, observe that every Motzkin path $\pi$ with $\tau_{2k}(\pi) = r$ can be constructed by the following procedure in exactly one way.
\begin{enumerate}
\item
Choose $d_0, \ldots, d_b$ and $l_0, \ldots, l_b$ such that (\ref{eqn:bigreven}) holds.
Construct a sequence of down and level steps which contains exactly $d_j$ down steps at height $k+j$ and $l_j$ level steps at height $k+j$ for $0 \le j \le b$ and which satisfies all of the following.
\begin{enumerate}
\item
The step immediately preceeding a step at height $j$ is either a down step at height $j+1$ or less or a level step at height $j$ or less.
\item
All steps after the last down step at height $j$ are at height $j-1$ or less.
\item
The sequence ends with a down step at height $k$.
\end{enumerate}
\item
If the first step is at height $k+j$, insert $j+1$ up steps before the first step.
Similarly, after each step except the last, insert enough up steps to reach the height of the next level or down step.
\item
After each down step at height $k$ except the last, insert an (possibly empty) upside-down Motzkin path of height at most $k-1$.
\item
Before the first step insert a path from height 0 to height $k-1$ which does not exceed height $k-1$.
\item
After the last step, insert a path from height $k-1$ to height 0 which does not exceed height $k-1$.
\end{enumerate}
Since the choice at each step is independent of the choices at the other steps, and since every sequence of choices results in a path of the type desired, the desired generating function is the product of the generating functions for each step.

To compute the generating function for step 1, suppose we have fixed $d_0, \ldots, d_b$ and $l_0, \ldots, l_b$;  then each of the resulting partial paths will have generating function $i^{\sum\limits_{j=0}^b d_j} x^{\sum\limits_{j=0}^b (d_j+l_j)}$.
To count these paths, we construct them from the top down.
That is, we first arrange the $l_b$ level steps at height $k+b$;  there is one way to do this.
We then place the $d_b$ down steps at height $k+b$ so that one of these steps occurs after all of the diagonal steps.
There are ${{d_b + l_b - 1} \choose {l_b}}$ ways to do this.
We then place the $l_{b-1}$ level steps at height $k+b-1$ so that none of these steps immediately follows a diagonal step at height $k+b$.
There are ${{d_b + l_{b-1}}\choose{l_{b-1}}}$ ways to do this.
Proceeding in this fashion, we find that the generating function for step 1 is equal to
\begin{equation}
\label{eqn:step1gf}
\sum \prod_{s=0}^b {{d_s + d_{s+1} + l_s - 1}\choose{d_{s+1} + l_s}} {{d_{s+1} + l_s}\choose{l_s}} i^{\sum\limits_{j=0}^b d_j} x^{\sum\limits_{j=0}^b (d_j + l_j)},
\end{equation}
where the sum on the left is over all sequences $d_0, \ldots, d_b$ and $l_0, \ldots, l_b$ of nonnegative integers which satisfy (\ref{eqn:bigreven}).
In the path obtained after step 2 there is exactly one up step for every down step so the generating function for step 2 is equal to
\begin{equation}
\label{eqn:step2gf}
(ix)^{\sum\limits_{j=0}^b d_j}.
\end{equation}
Using (\ref{eqn:Agf}) and (\ref{eqn:Auk}), we find that the generating function for step 3 is equal to
\begin{equation}
\label{eqn:step3gf}
\left( \frac{W_{k-1}(x)}{ix W_k(x)}\right)^{d_0-1}
\end{equation}
and the generating functions for steps 4 and 5 are both equal to
\begin{equation}
\label{eqn:step45gf}
\frac{1}{ix W_k(x)}.
\end{equation}
Taking the product of the quantities in (\ref{eqn:step1gf}), (\ref{eqn:step2gf}), (\ref{eqn:step3gf}) and the square of the quantity in (\ref{eqn:step45gf}), we obtain the desired generating function.
$\Box$

The proof of Theorem \ref{thm:biggfodd} is similar to the proof of Theorem \ref{thm:biggfeven}, using (\ref{eqn:Bgf}), (\ref{eqn:Cgf}), and (\ref{eqn:Buk}).

Theorems \ref{thm:biggfeven} and \ref{thm:biggfodd} have several interesting special cases;  we give two of them here.

\begin{corollary}
For all $k \ge 1$,
$$\sum_{\pi \in I(3412)} (-1)^{sign(\pi)} x^{|\pi|} = \frac{-1}{\left(W_k(x)\right)^2},$$
where the sum on the left is over all involutions in $I(3412)$ which contain exactly one subsequence of type $[2k]$.
\end{corollary}

\begin{corollary}
\label{cor:3211}
For all $k \ge 1$,
$$\sum_{\pi \in I(3412)} (-1)^{sign(\pi)} x^{|\pi|} = \frac{1}{i x \left( W_{k+1}(x) + W_k(x)\right)^2},$$
where the sum on the left is over all involutions in $I(3412)$ which contain exactly one subsequence of type $[2k+1]$.
\end{corollary}

\section{Involutions Which Avoid 3412 and Contain Another Pattern Exactly Once}

For any permutation $\pi$ we write
$$G_\pi^+(x) = \sum_{\sigma \in I(3412; \pi)} x^{|\sigma|},$$
and
$$G_\pi^-(x) = \sum_{\sigma \in I(3412; \pi)} (-1)^{sign(\sigma)} x^{|\sigma|}.$$
In this section we give recurrence relations for $G_\pi^+(x)$ and $G^-_\pi(x)$ when $\pi$ is direct sum indecomposable and we use these recurrence relations to give $G_\pi^+(x)$ and $G^-_\pi(x)$ in terms of Chebyshev polynomials for various $\pi$.
We begin with the case in which $\pi$ begins with $|\pi|$ and ends with 1.

\begin{proposition}
\label{prop:Gkpi1}
For any permutation $\pi$ we have
\begin{equation}
\label{eqn:Gpi+}
G_{1 \ominus \pi \ominus 1}^+(x) = x^2 G_\pi^+(x) \left( F_{1 \ominus \pi \ominus 1}^+(x)\right)^2
\end{equation}
and
\begin{equation}
\label{eqn:Gpi-}
G_{1 \ominus \pi \ominus 1}^-(x) = -x^2 G_\pi^-(x) \left( F_{1 \ominus \pi \ominus 1}^-(x) \right)^2.
\end{equation}
\end{proposition}
\begin{proof}
To prove (\ref{eqn:Gpi+}), first observe that the set $I(3412; 1 \ominus \pi \ominus 1)$ can be partitioned into two sets:  the set $A_1$ of those involutions which begin with 1 and the set $A_2$ of those involutions which do not begin with 1.

In view of Proposition \ref{prop:permmap}(i), the set $A_1$ contributes $x G^+_{1 \ominus \pi \ominus 1}(x)$ to the desired generating function.

In view of Proposition \ref{prop:permmap}(ii), every permutation in $A_2$ has the form $\sigma_1 * \sigma_2$.
Since $1 \ominus \pi \ominus 1$ is direct sum indecomposable, $A_2$ can be partitioned into two sets:  the set $B_1$ of involutions in which $1 \ominus \pi \ominus 1$ occurs in $1 \ominus \sigma_1 \ominus 1$ and the set $B_2$ of involutions in which $1 \ominus \pi \ominus 1$ occurs in $\sigma_2$.
The set $B_1$ contributes $x^2 G_\pi^+(x) F_{1 \ominus \pi \ominus 1}^+(x)$ to the desired generating function and the set $B_2$ contributes $x^2 G_{1 \ominus \pi \ominus 1}^+(x) F_\pi^+(x)$.

Add the contributions of $A_1$, $B_1$, and $B_2$ and solve the resulting equation for $G_{1 \ominus \pi \ominus 1}^+(x)$ to obtain
$$G_{1 \ominus \pi \ominus 1}^+(x) = \frac{x^2 G_\pi^+(x) F_{1 \ominus \pi \ominus 1}^+(x)}{1-x-x^2 F^+_\pi(x)}.$$
Now (\ref{eqn:Gpi+}) follows from (\ref{eqn:Fkpi1plus}).

The proof of (\ref{eqn:Gpi-}) is similar to the proof of (\ref{eqn:Gpi+}).
\end{proof}

Arguing as in the proof of Proposition \ref{prop:Gkpi1}, one can show that if $\pi$ does not end with 1 then (\ref{eqn:Gpi+}) and (\ref{eqn:Gpi-}) hold when $1 \ominus \pi \ominus 1$ is replaced with $1 \ominus \pi$.
Similarly, one can show that if $\pi$ does not begin with $|\pi|$ then (\ref{eqn:Gpi+}) and (\ref{eqn:Gpi-}) hold when $1 \ominus \pi \ominus 1$ is replaced with $\pi \ominus 1$.

To complete our analysis of $G_\pi^+(x)$ and $G_\pi^-(x)$ for the case in which $\pi$ is direct sum indecomposable, observe that there are four possibilities for the form of $\pi$, corresponding to whether or not $\pi$ begins with $|\pi|$ and whether or not $\pi$ ends with 1.
Proposition \ref{prop:Gkpi1} and its analogues address three of these possibilities.
To address the fourth, suppose $\pi$ is direct sum indecomposable, $\pi$ does not begin with $|\pi|$, and $\pi$ does not end with 1.
Then $\overline{\pi} = \pi$ and by \cite[Proposition 5.7]{Egge} every involution which avoids 3412 also avoids $\pi$.
It follows that $G^+_\pi(x) = G_\pi^-(x) = 0$.

For the remainder of this section we use Proposition \ref{prop:Gkpi1} and its analogues to compute $G_\pi^+(x)$ and $G^-_\pi(x)$ for various $\pi$.
In each case we express $G_\pi^+(x)$ and $G^-_\pi(x)$ in terms of Chebyshev polynomials of the second kind.
We begin with $G_{12}^+(x)$ and $G^-_{12}(x)$.

\begin{proposition}
We have
\begin{equation}
\label{eqn:G+12}
G^+_{12}(x) = \frac{x^2}{1-x^2}
\end{equation}
and
\begin{equation}
\label{eqn:G-12}
G^-_{12}(x) = \frac{x^2}{1+x^2}.
\end{equation}
\end{proposition}
\begin{proof}
First observe that if a permutation contains exactly one subsequence of type 12 then it has the form $n\ (n-1)\ldots i\ (i+1)\ \ldots 21$ for some $i$, $1 \le i \le n-1$.
Now observe that permutations of this form are involutions if and only if $n$ is even and $i = \frac{n}{2}$.
Now (\ref{eqn:G+12}) and (\ref{eqn:G-12}) follow.
\end{proof}

Using our expressions for $G_{12}^+(x)$ and $G^-_{12}(x)$, we obtain expressions for $G^+_{[k] \ominus 231}(x)$, $G^+_{[k] \ominus 12}(x)$, $G^-_{[k] \ominus 231}(x)$, and $G^-_{[k] \ominus 12}(x)$.

\begin{proposition}
\label{prop:Gk231312}
For all $k \ge 1$ we have
\begin{equation}
\label{eqn:G+k231}
G^+_{[k] \ominus 231}(x) = G^+_{[k] \ominus 12}(x) = \frac{(1-x)}{(1+x) V_{k+1}^2(x)}
\end{equation}
and
\begin{equation}
\label{eqn:G-k231}
G^-_{[k] \ominus 231}(x) = G^-_{[k] \ominus 12}(x) = \frac{x^2 (1+x^2)}{\left((1+x^2) W_k(x) - ix(1+x) W_{k-1}(x)\right)^2}.
\end{equation}
\end{proposition}
\begin{proof}
To see that $G^+_{[k] \ominus 231}(x) = G^+_{[k] \ominus 12}(x)$, first observe that if $\pi = [k-1] \ominus 12$ then $[k] \ominus 231 = 1 \ominus \pi \ominus 1$ and $[k] \ominus 12 = 1 \ominus \pi$.
Now use (\ref{eqn:Gpi+}), \cite[Theorem 6.3]{Egge}, and the analogue of (\ref{eqn:Gpi+}) in which $1 \ominus \pi \ominus 1$ is replaced with $1 \ominus \pi$ to find that
\begin{eqnarray*}
G^+_{1 \ominus \pi \ominus 1}(x) &=&  x^2 G_\pi^+(x) \left( F_{1 \ominus \pi \ominus 1}^+(x) \right)^2 \\
&=& x^2 G_\pi^+(x) \left( F_{1 \ominus \pi}^+(x) \right)^2 \\
&=& G^+_{1 \ominus \pi}(x),
\end{eqnarray*}
as desired.

To prove (\ref{eqn:G+k231}) we show that $G^+_{[k] \ominus 12}(x)$ is equal to the quantity on the right for $k \ge 2$.
Arguing by induction on $k$, first observe that when $k=0$ the result is immediate from (\ref{eqn:G+12}).
Now suppose $k \ge 1$ and the result holds for $k-1$.
Set $\pi = [k-1] \ominus 12$ in the analogue of (\ref{eqn:Gpi+}) in which $1 \ominus \pi \ominus 1$ is replaced with $1 \ominus \pi$, use \cite[Theorem 6.2]{Egge} to replace $F^+_{\pi}(x)$ with $\frac{V_{k-2}(x)}{x V_{k-1}(x)}$, and use induction to eliminate $G_\pi^+(x)$ and obtain the quantity on the right side of (\ref{eqn:G+k231}), as desired.

The proof of (\ref{eqn:G-k231}) is similar to the proof of (\ref{eqn:G+k231}).
\end{proof}

Next we compute $G^+_{213}(x)$, $G^+_{132}(x)$, $G^-_{213}(x)$, and $G^-_{132}(x)$.

\begin{proposition}
We have
\begin{equation}
\label{eqn:G213132+}
G^+_{213}(x) = G^+_{132}(x) = \frac{x^3}{1-x-x^2}
\end{equation}
and
\begin{equation}
\label{eqn:G213132-}
G^-_{213}(x) = G^-_{132}(x) = -\frac{x^3}{1-x+x^2}.
\end{equation}
\end{proposition}
\begin{proof}
First observe that 132 is the reverse complement of 213 and that this operation preserves parity, so $G^+_{213}(x) = G^+_{132}(x)$ and $G^-_{213}(x) = G^-_{132}(x)$.

We now show that $G^+_{132}(x)$ is equal to the quantity on the right side of (\ref{eqn:G213132+}).
To do this, first observe that the set $I(3412; 132)$ can be partitioned into two sets:  the set $A_1$ of those involutions which have the form $1 \ominus \sigma$ and the set $A_2$ of those involutions which have the form $\sigma_1 * \sigma_2$.

If $\pi \in A_1$, so that $\pi = 1 \ominus \sigma$, then $\sigma$ contains exactly one subsequence of type 21 and no subsequences of type 132.
Since $\sigma$ contains exactly one subsequence of type 21, it must have the form $1 \ldots (i+1)\ i \ldots n$ for some $i$, $1 \le i \le n-1$.
The only permutation of this form which avoids 132 is $213\ldots n$.
It follows that the set $A_1$ contributes $\frac{x^3}{1-x}$ to $G^+_{132}(x)$.

If $\pi \in A_2$, so that $\pi = \sigma_1 * \sigma_2$, then $\sigma_2$ contains no subsequences of type 21, so $\sigma_2 = 12 \ldots k$.
It follows that $\sigma_1$ contains exactly one subsequence of type 132.
Therefore the set $A_2$ contributes $\frac{x^2}{1-x} G^+_{132}(x)$ to $G^+_{132}(x)$.

Add the contributions of $A_1$ and $A_2$ and solve the resulting equation for $G^+_{132}(x)$ to obtain (\ref{eqn:G213132+}).

The proof of (\ref{eqn:G213132-}) is similar to the proof of (\ref{eqn:G213132+}).
\end{proof}

Using our expressions for $G^+_{213}(x)$, $G^+_{132}(x)$, $G^-_{213}(x)$, and $G^-_{132}(x)$, we obtain expressions for $G_{[k] \ominus 213}^+(x)$, $G_{[k] \ominus 132}^+(x)$, $G_{[k] \ominus 213}^-(x)$, and $G_{[k] \ominus 132}^-(x)$.

\begin{proposition}
For all $k \ge 0$ we have
\begin{displaymath}
G_{[k] \ominus 213}^+(x) = G_{[k] \ominus 132}^+(x) = \frac{1-x-x^2}{x \left( V_{k+2}(x) + V_{k+1}(x)\right)^2}
\end{displaymath}
and
\begin{displaymath}
G_{[k] \ominus 213}^-(x) = G_{[k] \ominus 132}^-(x) = \frac{x(1-x+x^2)}{\left(W_{k+1}(x) - ix W_k(x)\right)^2}.
\end{displaymath}
\end{proposition}
\begin{proof}
This is similar to the second half of the proof of Proposition \ref{prop:Gk231312}, using (\ref{eqn:G213132+}) and (\ref{eqn:G213132-}).
\end{proof}

Next we compute $G^+_{123}(x)$ and $G^-_{123}(x)$.

\begin{proposition}
We have
\begin{equation}
\label{eqn:G123+}
G^+_{123}(x) = \frac{x^3(1+x^2)}{(1-x^2)^2}
\end{equation}
and
\begin{equation}
\label{eqn:G123-}
G^-_{123}(x) = \frac{x^3(1-x^2)}{(1+x^2)^2}.
\end{equation}
\end{proposition}
\begin{proof}
To prove (\ref{eqn:G123+}), first observe that the set $I(3412; 123)$ can be partitioned into two sets:  the set $A_1$ of those involutions which have the form $1 \ominus \sigma$ and the set $A_2$ of those involutions which have the form $\sigma_1 * \sigma_2$.

If $\pi \in A_1$, so that $\pi = 1 \ominus \sigma$, then $\sigma$ contains exactly one subsequence of type 12.
It follows from (\ref{eqn:G+12}) that $A_1$ contributes $\frac{x^2}{1-x^2}$ to $G^+_{123}(x)$.

If $\pi \in A_2$, so that $\pi = \sigma_1 * \sigma_2$, then we must have $|\sigma_2| \le 1$.
If $|\sigma_2| = 0$ then $\sigma_1$ contains exactly one subsequence of type 123.
If $|\sigma_2| = 1$ then $\sigma_1$ contains exactly one subsequence of type 12.
It follows from (\ref{eqn:G+12}) that $A_2$ contributes $x^2 G_{123}^+(x) + \frac{x^5}{1-x^2}$ to $G^+_{123}(x)$.

Add the contributions of $A_1$ and $A_2$ and solve the resulting equation for $G_{123}^+(x)$ to obtain (\ref{eqn:G123+}).

The proof of (\ref{eqn:G123-}) is similar to the proof of (\ref{eqn:G123+}).
\end{proof}

Using our expressions for $G_{123}^+(x)$ and $G_{123}^-(x)$, we obtain expressions for $G^+_{[k] \ominus 123}(x)$ and $G^-_{[k] \ominus 123}(x)$.

\begin{proposition}
For all $k \ge 1$ we have
\begin{displaymath}
G^+_{[k] \ominus 123}(x) = \frac{x (1+x^2)(1-2x+2x^3-x^4)^2}{(1-x^2)^2 \left( (1-x+x^3) V_{k+1}(x) + x (x-1) V_k(x)\right)^2}
\end{displaymath}
and
\begin{displaymath}
G^-_{[k] \ominus 123}(x) = \frac{x^3 (1-x^2)(1+x^2)^4}{\left((1+x^2)^3 W_k(x) - ix(1 + 3x^2+2x^4+x^5) W_{k-1}(x)\right)^2}.
\end{displaymath}
\end{proposition}
\begin{proof}
This is similar to the second half of the proof of Proposition \ref{prop:Gk231312}, using (\ref{eqn:G123+}) and (\ref{eqn:G123-}).
\end{proof}

\section{Involutions Which Contain 3412 Exactly Once and Avoid Another Pattern}

In this section we study involutions which contain exactly one subsequence of type 3412.
We begin by describing the recursive structure of the set of such involutions.
We then use this recursive structure to find the generating functions
$$P_\pi^+(x) = \sum_{\sigma \in I(\pi; 3412)} x^{|\sigma|}$$
and
$$P_\pi^-(x) = \sum_{\sigma \in I(\pi; 3412)} (-1)^{sign(\sigma)} x^{|\sigma|}$$
for various permutations $\pi$.
We start with the notion of a crossing in an involution.

\begin{definition}
Let $\pi$ denote an involution.
A {\em crossing} in $\pi$ is a sequence $i < j < k < l$ such that $\pi(i) = k$ and $\pi(j) = l$.
Each crossing has a corresponding subsequence in $\pi$, given by $\pi(i) \pi(j) \pi(k) \pi(l)$.
\end{definition}

As we show next, involutions which avoid 3412 are exactly those involutions which have no crossings.

\begin{proposition}
\label{prop:3412crossings}
Suppose $\pi$ is an involution.
Then $\pi$ avoids 3412 if and only if $\pi$ has no crossings.
\end{proposition}
\begin{proof}
($\Longrightarrow$)
Suppose $\pi \in I(3412)$.
Observe that if $abcd$ is a crossing then its corresponding subsequence has type 3412.
But $\pi$ avoids 3412, so $\pi$ has no crossings.

($\Longleftarrow$)
Suppose $\pi$ is an involution and $\pi(a)$ $\pi(b)$ $\pi(c)$ $\pi(d)$ is a subsequence of type 3412.
Furthermore, suppose $\pi(a) \le a$.
Since $\pi(c) < \pi(a) \le a$ and $\pi(d) < \pi(a) \le a$, we find that $\pi(c) < c$ and $\pi(d) < d$.
It follows that $\pi(c) \pi(d) c d$ is a crossing, since $\pi(c) < \pi(d)$.
Therefore $\pi(a) > a$, and by a similar argument we find $\pi(b) > b$.
Now it follows that $a b \pi(a) \pi(b)$ is a crossing, since $\pi(a) < \pi(b)$.
\end{proof}

We now consider crossings in involutions which contain exactly one subsequence of type 3412.

\begin{corollary}
\label{cor:thecrossing}
Suppose $\pi$ is an involution with exactly one subsequence of type 3412, given by $\pi(a) \pi(b) \pi(c) \pi(d)$.
Then $\pi$ has exactly one crossing, which is $abcd$.
\end{corollary}
\begin{proof}
By Proposition \ref{prop:3412crossings}, the involution $\pi$ has at least one crossing.
Moreover, each crossing corresponds to a subsequence of type 3412, so $\pi$ has at most one crossing.
It follows that $\pi$ has exactly one crossing.
In addition, the subsequence corresponding to this crossing has type 3412, so the subsequence must be $\pi(a) \pi(b) \pi(c) \pi(d)$.
It follows that the crossing is $abcd$.
\end{proof}

We now describe the recursive structure of $I(;3412)$.

\begin{theorem}
\label{thm:Jncases}
Fix $\pi \in I_n(;3412)$.
Then exactly one of the following holds.
\renewcommand\labelenumi{{\upshape (\roman{enumi}) }}
\begin{enumerate}
\item
$\pi$ begins with 1.
\item
$\pi$ does not begin with 1, 1 is not an element of the subsequence of $\pi$ of type 3412, and this subsequence lies between $\pi(1)$ and 1.
\item
$\pi$ does not begin with 1, 1 is not an element of the subsequence of $\pi$ of type 3412, and this subsequence lies entirely to the right of 1.
\item
1 is an element of the subsequence of $\pi$ of type 3412.
\end{enumerate}
\end{theorem}
\begin{proof}
Observe that it is sufficient to prove both of the following.
\begin{enumerate}
\item[(A)]
If $\pi$ does not begin with 1 and 1 is not an element of the subsequence of type 3412 then this subsequence lies between $\pi(1)$ and 1 or it lies entirely to the right of 1.
\item[(B)]
If 1 is an element of the subsequence of $\pi$ of type 3412 then $\pi$ does not begin with 1.
\end{enumerate}

To prove statement (A), suppose by way of contradiction that $\pi(a) \pi(b) \pi(c) \pi(d)$ is the subsequence of type 3412 and 1 lies between $\pi(a)$ and $\pi(d)$.
By Corollary \ref{cor:thecrossing} we have $\pi(a) = c$ and $\pi(b) = d$.
It follows that at least one of $\pi(1) a 1 c$ and $\pi(1) b 1 d$ is a crossing in $\pi$.
But $abcd$ is the unique crossing in $\pi$, so $c = \pi(a) = 1$, which contradicts the fact that 1 is not an element of the subsequence of $\pi$ of type 3412.

To prove statement (B), observe that if 1 is an element of the subsequence of type 3412 in $\pi$ then it is the third element of this subsequence, so it cannot be the first element of $\pi$.
\end{proof}

\begin{definition}
For any permutations $\pi_1$, $\pi_2$, and $\pi_3$ we write $\pi_1 \otimes \pi_2 \otimes \pi_3$ to denote the permutation in $S_{|\pi_1|+|\pi_2|+|\pi_3|+4}$ which is given by
$$(\pi_1 \otimes \pi_2 \otimes \pi_3)(i) =
\cases{
|\pi_1|+|\pi_2|+3 & if $i = 1$ \cr
\pi_1(i-1) + 1 & if $2 \le i \le |\pi_1|+1$ \cr
|\pi_1| + |\pi_2| + |\pi_3|+4 & if $i = |\pi_1| + 2$ \cr
\pi_2(i - |\pi_1| - 2) + |\pi_1| + 2 & if $3 \le i - |\pi_1| \le |\pi_2| + 2$ \cr
1 & if $i = |\pi_1| + |\pi_2| +3$ \cr
\pi_3(i - |\pi_1| - |\pi_2| - 3) + |\pi_1| + |\pi_2| + 3 & if $4 \le i - |\pi_1| - |\pi_2| \le + |\pi_3| + 3$ \cr
1 + |\pi_1| & if $i = |\pi_1| + |\pi_2| + |\pi_3| + 4$ \cr
}.$$
\end{definition}

\begin{theorem}
\label{thm:Jndecomp}
\renewcommand\labelenumi{{\upshape (\roman{enumi}) }}
\begin{enumerate}
\item
For all $n \ge 1$, the map
$$
\begin{array}{rcl}
I_{n-1}(;3412) & \longrightarrow & I_n(;3412) \\
\pi &\mapsto& 1 \oplus \pi \\
\end{array}
$$
is a bijection between $I_{n-1}(;3412)$ and the set of involutions in $I_n(;3412)$ which begin with 1.
\item
For all $n \ge 1$ and all $j$ such that $2 \le j \le n$, the map
$$
\begin{array}{rcl}
I_{j-2}(;3412) \times I_{n-j}(3412) &\longrightarrow& I_n(;3412) \\
(\pi_1, \pi_2) &\mapsto& \pi_1 * \pi_2 \\
\end{array}
$$
is a bijection between $I_{j-2}(;3412) \times I_{n-j}(3412)$ and the set of involutions in $I_n(;3412)$ which do not begin with 1, in which 1 is not an element of the subsequence of type 3412, and in which this subsequence lies between $\pi(1)$ and 1.
\item
For all $n \ge 1$ and all $j$ such that $2 \le j \le n$, the map
$$
\begin{array}{rcl}
I_{j-2}(3412) \times I_{n-j}(;3412) &\longrightarrow& I_n(;3412) \\
(\pi_1, \pi_2) &\mapsto& \pi_1 * \pi_2 \\
\end{array}
$$
is a bijection between $I_{j-2}(3412) \times I_{n-j}(3412)$ and the set of involutions in $I_n(;3412)$ which do not begin with 1, in which 1 is not an element of the subsequence of type 3412, and in which this subsequence lies entirely to the right of 1.
\item
For each $k \ge 0$, set $[k] = k\ k-1\ \ldots 2 1$.
For all $n \ge 1$ and all $j,k,m \ge 0$ such that $0 \le j+k+m \le n-4$, the map
$$
\begin{array}{ccc}
I_j(3412) \times I_k(3412) \times I_m(3412) &\longrightarrow& I_n(;3412)
\end{array}
$$
$$
\begin{array}{ccc}
(\pi_1, \pi_2, \pi_3) &\mapsto& (\pi_1 \otimes [n-j-k-m-4] \otimes \pi_2) \oplus \pi_3
\end{array}
$$
is a bijection between $I_j(3412) \times I_k(3412) \times I_m(3412)$ and the set of involutions in $I_n(;3412)$ in which 1 is an element of the subsequence of type 3412.
\end{enumerate}
\end{theorem}
\begin{proof}
(i)
This is immediate from the fact that if $\pi$ begins with 1 then 1 cannot be an element of a subsequence of type 3412 in $\pi$.

(ii)
Observe that if the subsequence of type 3412 lies between $\pi(1)$ and 1 then the elements of $\pi$ to the left of 1 are $2,3,\ldots,\pi(1)$.
Now (ii) follows.

(iii)
This is similar to the proof of (ii).

(iv)
It is routine to verify that every permutation of the form given contains exactly one subsequence of type 3412, and that 1 is an element of this subsequence.

To show that the given map is one-to-one and onto, suppose $\pi \in I_n(;3412)$ and that 1 is an element of the subsequence of type 3412 in $\pi$.
Let $\pi(a) \pi(b) \pi(c) \pi(d)$ denote this subsequence.
Then there exist unique sequences $\sigma_0, \sigma_1, \sigma_2, \sigma_3, \sigma_4$ which avoid 3412 such that $\pi = \sigma_0 a \sigma_1 b \sigma_2 c \sigma_3 d \sigma_4$.
Since 1 is an element of $\pi(a) \pi(b) \pi(c) \pi(d)$ and this sequence has type 3412, we have $\pi(c) = 1$.
By Corollary \ref{cor:thecrossing} we have $\pi(1) = a$, so $\sigma_0$ is empty.

Now set $j = |\sigma_1|$, $k = |\sigma_3|$, and $m = |\sigma_4|$, so that $|\sigma_2| = n-j-k-m-4$.
Since $abcd$ is the only crossing in $\pi$, the entries of $\sigma_1$ are $2, \ldots, j+1$, the entries of $\sigma_2$ are $j+3, \ldots, n-k-m-2$, the entries of $\sigma_3$ are $n-k-m, \ldots, n-m-1$, and the entries of $\sigma_4$ are $n-m+1, \ldots, n$.
Moreover, if $\sigma_2$ has a subsequence $\pi(e) \pi(f)$ of type 12 then $\pi(a) \pi(b) \pi(e) \pi(f)$ will have type 3412.
It follows that $\sigma_2$ has type $[n-j-k-m-4]$.
Combining these observations, we find there exist unique permutations $\pi_1, \pi_2, \pi_3$ such that $|\pi_1| = j$, $|\pi_2| = k$, $|\pi_3| = m$, and $\pi = (\pi_1 \otimes [n-j-k-m-4] \otimes \pi_2) \oplus \pi_3$.
It follows that the given map is a bijection, as desired.
\end{proof}

Using Theorem \ref{thm:Jndecomp}, we find the generating function for those involutions which contain exactly one subsequence of type 3412.

\begin{proposition}
We have
\begin{equation}
\label{eqn:P+empty}
P_{\emptyset}^+(x) = \frac{2x-1}{2x^2(1-x)} + \frac{1-2x-2x^2}{2x^2 \sqrt{1-2x-3x^2}}
\end{equation}
and
\begin{equation}
\label{eqn:P-empty}
P^-_{\emptyset}(x) = \frac{(x+1)(2x^2-2x+1)}{2x^2 (1+x^2)} + \frac{(x^2-1)(4x^2-2x+1)}{2x^2(1+x^2)\sqrt{5x^2-2x+1}}.
\end{equation}
\end{proposition}
\begin{proof}
To prove (\ref{eqn:P+empty}), first use Theorems \ref{thm:Jncases} and \ref{thm:Jndecomp} to find that
$$P_{\emptyset}^+(x) = x P_\emptyset^+(x) + 2 x^2 P_\emptyset^+(x) F_\emptyset^+(x) + \frac{x^4 (F_\emptyset^+(x))^3}{1-x}.$$
Now use the fact that
$$F_\emptyset^+(x) = \frac{1-x-\sqrt{1-2x-3x^2}}{2x^2}$$
to eliminate $F_\emptyset^+(x)$ and solve the resulting equation for $P^+_\emptyset(x)$ to obtain (\ref{eqn:P+empty}).

The proof of (\ref{eqn:P-empty}) is similar to the proof of (\ref{eqn:P+empty}), using
$$P^-_\emptyset(x) = x P_\emptyset^-(x) - 2x^2 P_\emptyset^-(x) F_\emptyset^-(x) + \frac{x^4 (1+x) (F_\emptyset^-(x))^3}{1+x^2}$$
and
$$F_\emptyset^-(x) = \frac{x-1+\sqrt{5x^2-2x+1}}{2x^2}.$$
\end{proof}

For the rest of this section we use Theorem \ref{thm:Jndecomp} to find $P^+_\pi(x)$ and $P^-_\pi(x)$ for various $\pi$.
We begin with the case in which $\pi = [k]$.

\begin{proposition}
\label{prop:P+k21}
For all $k \ge 1$ we have
\begin{equation}
\label{eqn:P+2k}
P^+_{[2k]}(x) = \frac{\sum_{j=0}^{k-2}(1-x^{2j+2}) V_j^2(x)}{(1-x) V_k^2(x)}
\end{equation}
and
\begin{equation}
\label{eqn:P+2k+1}
P^+_{[2k+1]}(x) = \frac{\sum_{j=0}^{k-1} (1-x^{2j+1}) \left( V_j(x) + V_{j-1}(x)\right)^2}{(1-x)\left(V_{k+1}(x) + V_k(x)\right)^2}.
\end{equation}
\end{proposition}
\begin{proof}
To begin, use Theorems \ref{thm:Jncases} and \ref{thm:Jndecomp} to find that
\begin{eqnarray*}
\lefteqn{P_{[k]}^+(x) =} & & \\
& & x P^+_{[k]}(x) + x^2 P_{[k-2]}(x) F_{[k]}^+(x) + x^2 F_{[k-2]}^+(x) P_{[k]}^+(x) \\
& & + x^4 (1 + x + \cdots +x^{k-3}) \left( F_{[k-2]}^+(x)\right)^2 F_{[k]}^+(x).
\end{eqnarray*}
Now solve this equation for $P_{[k]}^+(x)$, use \cite[Theorem 6.1]{Egge} and argue by induction on $k$ to obtain (\ref{eqn:P+2k}) and (\ref{eqn:P+2k+1}).
\end{proof}

Proposition \ref{prop:P+k21} enables us to find $|I_n([k];3412)|$ for various $k$.
We obtain simple formulas for this quantity when $k = 3$ or $k = 4$.

\begin{corollary}
For all $n \ge 4$,
$$|I_n(321;3412)| = \frac{n-3}{5}F_{n-1} + \frac{n-1}{5} F_{n-3}$$
and
$$|I_n(4321;3412)| = 2^{n-5} (3n-10).$$
Here $F_n$ is the $n$th Fibonacci number, which is defined by $F_0 = 0$, $F_1 = 1$, and $F_n = F_{n-1}+F_{n-2}$ for $n \ge 2$.
\end{corollary}
\begin{proof}
Set $k = 1$ in (\ref{eqn:P+2k+1}) and set $k = 2$ in (\ref{eqn:P+2k}).
\end{proof}

Arguing as in the proof of Proposition \ref{prop:P+k21}, we also obtain $P^-_{[k]}(x)$.

\begin{proposition}
\label{prop:P-k21}
For all $k \ge 1$ we have
\begin{equation}
\label{eqn:P-2k}
P^-_{[2k]}(x) = \frac{\sum_{j=0}^{k-1} v_{2j}(x) W_j^2(x)}{W_k^2(x)}
\end{equation}
and
\begin{equation}
\label{eqn:P-2k+1}
P^-_{[2k+1]}(x) = \frac{\sum_{j=0}^{k-1} v_{2j}(x) \left( W_j(x) - i W_{j-1}(x) \right)^2}{\left( W_{k+1}(x) - i W_k(x)\right)^2}.
\end{equation}
Here $v_k(x) = \sum_{j=0}^k (-1)^{{{j}\choose{2}}} x^j$ for all $k \ge 0$.
\end{proposition}

Propositions \ref{prop:P+k21} and \ref{prop:P-k21} include $P^+_{321}(x)$ and $P^-_{321}(x)$ as special cases.
We now use Theorem \ref{thm:Jndecomp} to compute $P^+_\pi(x)$ and $P^-_\pi(x)$ for all remaining permutations of length 3.

\begin{proposition}
\label{prop:P+S3}
We have

\renewcommand{\baselinestretch}{2}
\small
\normalsize

\bigskip

\begin{center}
\begin{tabular}{|c|c|c|}
\hline
$\pi$ & $P_\pi^+(x)$ & $P_\pi^-(x)$ \\
\hline
$123$ & ${\displaystyle \frac{x^4(1+x+x^2)}{(1+x)(1-x)^3}}$ & ${\displaystyle \frac{x^4(1+3x-3x^2+6x^3+3x^4+3x^5-x^6)}{(1+x^2)(1+x)^3(1-x)^3}}$ \\
\hline
$132, 213$ & ${\displaystyle \frac{x^4}{(1-x)(1-x-x^2)}}$ & ${\displaystyle \frac{x^4(1+x)}{(1-x+x^2)(1+x^2)}}$ \\
\hline
$231, 312$ & $0$ & $0$ \\
\hline
\end{tabular}
\end{center}

\renewcommand{\baselinestretch}{1}
\small
\normalsize

\end{proposition}
\begin{proof}
The last line of the table is immediate, since a permutation which contains a subsequence of type 3412 must also contain subsequences of type 231 and 312.

The fact that $P_{132}^+(x) = P_{213}^+(x)$ is immediate, since $213$ is the reverse complement of $132$.

To obtain $P_{132}^+(x)$, use Theorem \ref{thm:Jndecomp} to find that
$$P_{132}^+(x) = x^2 \left( \frac{1}{1-x} \right) P_{132}^+(x) + x^4 \cdot 1 \cdot \frac{1}{1-x} \cdot \frac{1}{1-x}.$$
Solve this equation for $P_{132}^+(x)$ to obtain the desired result.

The rest of the table can be obtained in a similar fashion.
\end{proof}

Building on Proposition \ref{prop:P+S3}, we now compute $P_{[k] \ominus \pi}^+(x)$ and $P_{[k] \ominus \pi}^-(x)$ for various $\pi$ of length 2 and 3.
We begin with $\pi = 12$.

\begin{proposition}
For all $k \ge 1$ we have
\begin{equation}
\label{eqn:Pk312+}
P_{[k] \ominus 12}^+(x) = \frac{\sum_{j=1}^{k-1} (1-x^j) V_j^2(x)}{(1-x) V_{k+1}^2(x)}
\end{equation}
and
\begin{equation}
\label{eqn:Pk312-}
P_{[k] \ominus 12}^-(x) = \frac{\sum_{j=1}^{k-1} v_{j-1}(x) \left( W_{j-2}(x) - ix f_0(x) W_{j-3}(x)\right)^2}{\left( W_{k-1}(x) - ix f_0(x) W_{k-2}(x)\right)^2}.
\end{equation}
Here $v_k(x) = \sum_{j=0}^k (-1)^{{{j}\choose{2}}} x^j$ for all $k \ge 0$ and $f_0(x) = \frac{1+x^2}{1-x+2x^2}$.
\end{proposition}
\begin{proof}
To prove (\ref{eqn:Pk312+}) we argue by induction on $k$.
The case $k = 1$ is immediate from Proposition \ref{prop:P+S3}, so suppose $k \ge 2$ and the result holds for $k-1$.
Use Theorems \ref{thm:Jncases} and \ref{thm:Jndecomp} to find that
\begin{eqnarray*}
\lefteqn{P_{[k] \ominus 12}^+(x) =} \\
& & x P^+_{[k] \ominus 12}(x) + x^2 P_{[k-1] \ominus 12}^+(x) + x^2 F_{[k-1] \ominus 12}^+(x) P_{[k] \ominus 12}^+(x) \\
& & + x^4 \frac{1-x^{k-3}}{1-x} \left(F_{[k-1] \ominus 12}^+(x)\right)^2 F_{[k] \ominus 12}^+(x).\\
\end{eqnarray*}
Solve this equation for $P_{[k] \ominus 12}^+(x)$ and use induction and the fact \cite[Theorem 6.2]{Egge} that $F_{[k] \ominus 12}^+(x) = \frac{V_{k-2}(x)}{x V_{k-1}(x)}$ to obtain (\ref{eqn:Pk312+}).

The proof of (\ref{eqn:Pk312-}) is similar to the proof of (\ref{eqn:Pk312+}).
\end{proof}

We now turn our attention to $P^+_{[k] \ominus 132}(x)$ and $P^+_{[k] \ominus 213}(x)$.

\begin{proposition}
For all $k \ge 1$ we have
\begin{equation}
\label{eqn:Pk132+}
P^+_{[k] \ominus 132}(x) = P_{[k] \ominus 213}^+(x) = \frac{(1-x^3) \left( V_2(x) + V_1(x) \right) + \sum_{j=2}^k \left( V_j(x) + V_{j-1}(x) \right)^2}{(1-x) \left( V_{k+2}(x) + V_{k+1}(x)\right)^2}.
\end{equation}
\end{proposition}
\begin{proof}
It is immediate that $P^+_{[k] \ominus 132}(x) = P_{[k] \ominus 213}^+(x)$, since $[k] \ominus 132$ is the reverse complement of $[k] \ominus 213$.
To prove that $P^+_{[k] \ominus 132}(x)$ is equal to the expression on the right side of (\ref{eqn:Pk132+}) we argue by induction on $k$.

To prove the result when $k = 1$, use Theorems \ref{thm:Jncases} and \ref{thm:Jndecomp} to find that
\begin{eqnarray*}
\lefteqn{P^+_{4132}(x) =} \\
& & x P_{4132}^+(x) + x^2 P_{132}^+(x) F_{4132}^+(x) + x^2 F_{132}^+(x) P_{4132}^+(x) + \\
 & & x^4 \cdot (F_{132}^+(x) - 1) \cdot 1 \cdot 1 \cdot F_{4132}^+(x) + x^4 \cdot 1 \cdot \frac{1}{1-x} \cdot 1 \cdot F_{4132}^+(x).\\
\end{eqnarray*}
Solve this equation for $P_{4132}^+(x)$ and use Proposition \ref{prop:P+S3} and the fact \cite[Theorem 6.4]{Egge} that $F_{132}^+(x) = \frac{1}{1-x-x^2}$ and $F_{4132}^+(x) = \frac{1-x-x^2}{1-2x-x^2+x^3}$ to find that (\ref{eqn:Pk132+}) holds when $k = 1$.

Now suppose $k \ge 2$ and (\ref{eqn:Pk132+}) holds for $k - 1$.
Use Theorems \ref{thm:Jncases} and \ref{thm:Jndecomp} to find that
\begin{eqnarray*}
\lefteqn{P^+_{[k] \ominus 132}(x) =} \\
& & x P_{[k] \ominus 132}^+(x) + x^2 P^+_{[k-1] \ominus 132}(x) F_{[k] \ominus 132}^+(x) + x^2 F_{[k-1] \ominus 132}^+(x) P_{[k] \ominus 132}^+(x) + \\
& & x^4 F^+_{[k-1] \ominus 132}(x) \cdot \frac{1}{1-x} \cdot F_{[k-1] \ominus 132}^+(x) \cdot F_{[k] \ominus 132}^+(x) \\
\end{eqnarray*}
Solve this equation for $P^+_{[k] \ominus 132}(x)$ and use induction and \cite[Theorem 6.4]{Egge} to obtain (\ref{eqn:Pk132+}).
\end{proof}

In order to find $P^+_{[k] \ominus 231}(x)$, we first consider a more general situation.

\begin{proposition}
\label{prop:P+k1}
For any permutation $\pi$ we have
\begin{equation}
\label{eqn:P+kpi1}
P^+_{1 \ominus \pi \ominus 1}(x) = x^2 \left(F_{1 \ominus \pi \ominus 1}^+(x)\right)^2 \left( P_\pi^+(x) + \frac{x^2}{1-x} \left(F_\pi^+(x)\right)^2 \right)
\end{equation}
and
\begin{equation}
\label{eqn:P-kpi1}
P^-_{1 \ominus \pi \ominus 1}(x) = -x^2 \left( F_{1 \ominus \pi \ominus 1}^-(x) \right)^2 \left( P_\pi^-(x) - \frac{x^2(1+x)}{1+x^2} \left(F_\pi^-(x)\right)^2 \right).
\end{equation}
\end{proposition}
\begin{proof}
To prove (\ref{eqn:P+kpi1}), use Theorems \ref{thm:Jncases} and \ref{thm:Jndecomp} to find that
$$\left(1-x-x^2 F_\pi^+(x)\right) P_{1 \ominus \pi \ominus 1}^+(x) = x^2 F_\pi^+(x) P_\pi^+(x) + \frac{x^4}{1-x} \left(F_\pi^+(x)\right)^2 F^+_{1 \ominus \pi \ominus 1}(x).$$
Now (\ref{eqn:P+kpi1}) follows from (\ref{eqn:Fkpi1plus}).

The proof of (\ref{eqn:P-kpi1}) is similar to the proof of (\ref{eqn:P+kpi1}).
\end{proof}

Arguing as in the proof of Proposition \ref{prop:P+k1}, one can show that if $\pi$ does not end with 1 then (\ref{eqn:P+kpi1}) and (\ref{eqn:P-kpi1}) hold when $1 \ominus \pi \ominus 1$ is replaced with $1 \ominus \pi$.
Similarly, one can also show that if $\pi$ does not begin with $|\pi|$ then (\ref{eqn:P+kpi1}) and (\ref{eqn:P-kpi1}) hold when $1 \ominus \pi \ominus 1$ is replaced with $\pi \ominus 1$.

\begin{proposition}
For all $k \ge 2$ we have
$$P^+_{[k] \ominus 231}(x) = \frac{V^2_{k-1}(x) + \sum_{j=1}^{k-2} (1-x^j) V_j^2(x)}{(1-x)V_{k+1}^2(x)}.$$
\end{proposition}
\begin{proof}
Set $\pi = [k-1] \ominus 12$ in (\ref{eqn:P+kpi1}) and use (\ref{eqn:Pk312+}) and \cite[Theorem 6.2]{Egge} to simplify the result.
\end{proof}

In view of Conjecture \ref{conj:minussymmetry} and \cite[Conjecture 6.9]{Egge}, one might conjecture that $P^+_{[l_1,\ldots,l_m]}(x)$ and $P^-_{[l_1,\ldots,l_m]}(x)$ are symmetric in $l_1, \ldots, l_m$.
However, it is not difficult to verify (using Maple) that this fails for $[l_1,l_2,l_3] = [2,1,1]$ and $[l_1,l_2,l_3] = [1,2,1]$.

\section{Involutions Which Contain 3412 and Another Pattern Exactly Once}

For any permutation $\pi$ and any $n \ge 0$, let $I_n(;3412,\pi)$ (resp. $I(;3412,\pi)$) denote the set of involutions of length $n$ (resp. of any length) which contain exactly one subsequence of type 3412 and exactly one subsequence of type $\pi$.
In this section we begin a study of the generating functions
$$Q^+_\pi(x) = \sum_{\sigma \in I(;3412,\pi)} x^{|\sigma|}$$
and
$$Q^-_\pi(x) = \sum_{\sigma \in I(;3412,\pi)} (-1)^{sign(\sigma)} x^{|\sigma|}.$$
It seems possible to obtain results similar to the results of previous sections, but the expressions involved become cumbersome rather quickly.
For this reason we restrict our attention to $Q^+_{[k]}(x)$.

\begin{theorem}
\label{thm:Q+k}
For all $k \ge 2$ we have
\begin{equation}
\label{eqn:Q+2k}
Q^+_{[2k]}(x) = \frac{2x}{(1-x) V_k^2(x)} \sum_{i=2}^k \frac{\left(1 - x^{2(i-1)}\right) V_{i-2}(x) V_i(x) + \sum_{j=0}^{i-2} \left( 1 - x^{2(j+1)} \right) V_j^2(x)}{V_{i-1}(x) V_i(x)}
\end{equation}
and
\begin{displaymath}
Q^+_{[2k+1]}(x) = \frac{2}{(1-x) Y^2_{k+1}(x)} \sum_{i=2}^k \frac{\left(1-x^{2i+1}\right) Y_i(x) Y_{i+2}(x) + \sum_{j=0}^i \left( 1-x^{2j+1}\right) Y_j^2(x)}{Y_{i+1}(x) Y_{i+2}(x)},
\end{displaymath}
where $Y_k(x) = V_k(x) + V_{k-1}(x)$ for all $k$.
\end{theorem}
\begin{proof}
In view of Theorems \ref{thm:Jncases} and \ref{thm:Jndecomp} we have
\begin{eqnarray*}
\lefteqn{Q_{[k]}^+(x) =} & & \\
& & x Q^+_{[k]}(x) + x^2 Q^+_{[k-2]}(x) F_{[k]}^+(x) + x^2 P_{[k-2]}^+(x) G_{[k]}^+(x) + x^2 F_{[k-2]}^+(x) Q_{[k]}^+(x) \\[2ex]
& & + x^2 G^+_{[k-2]}(x) P_{[k]}^+(x) + x^4 \left( \frac{1-x^{2k-2}}{1-x}\right) \left( 2 G_{[k-2]}^+(x) F_{[k-2]}^+(x) F_{[k]}^+(x) + G_{[k]}^+(x) \left(F_{[k-2]}^+(x)\right)^2\right). \\
\end{eqnarray*}
Now the result follows by induction on $k$, using the fact that $Q_{[2]}^+(x) = 0$ along with \cite[Theorems 4.2, 4.3, 6.1]{Egge}, (\ref{eqn:P+2k}), and (\ref{eqn:P+2k+1}).
\end{proof}

Theorem \ref{thm:Q+k} enables us to find $|I_n(;3412,[k])|$ for various $k$.
We obtain a simple formula for this quantity when $k = 4$.

\begin{corollary}
For all $n \ge 6$ we have
$$|I_n(;3412,4321)| = (3n^2-23n+38)2^{n-8}.$$
\end{corollary}
\begin{proof}
Set $k = 2$ in (\ref{eqn:Q+2k}).
\end{proof}

\section{Directions for Future Research}

\begin{enumerate}
\item
For any $m \ge 1$ and any permutation $\pi \in I(3412)$, let
$$\omega_m(\pi) = \sum_{k=1}^\infty (-m)^{k-1} \tau_k (\pi).$$
Recall that $\omega_1(\pi)$ is both the number of left-to-right maxima in $\pi$ and the number of right-to-left minima in $\pi$, and that $\omega_2(\pi)$ is the number of fixed points in $\pi$.
Find a combinatorial interpretation of $\omega_m(\pi)$ for $m \ge 3$.

\item
In view of (\ref{eqn:F2kminus}), (\ref{eqn:Fk312}), and \cite[Theorems 6.1 and 6.2]{Egge}, the number of even (resp. odd) permutations in $I_n(3412,[2k])$ is equal to the number of even (resp. odd) permutations in $I_n(3412, [k-1] \ominus 231)$ for all $n\ge 0$.
Similarly, in view of (\ref{eqn:F2k-1minus}), (\ref{eqn:F213132}), and \cite[Theorems 6.1 and 6.4]{Egge}, the number of even (resp. odd) permutations in $I_n(3412,[2k+1])$ is equal to the number of even (resp. odd) permutations in $I_n(3412,[k-1] \ominus 132)$ for all $n\ge 0$.
Give combinatorial proofs of these facts.

\item
Prove Conjecture \ref{conj:minussymmetry}, which says that $F^-_{[l_1,\ldots,l_m]}(x)$ is symmetric in $l_1,\ldots, l_m$.
Alternatively, prove Conjecture \ref{conj:evenoddsymmetry}, which says that the generating functions for the even involutions in $I(3412, [l_1, \ldots, l_m])$ and the odd involutions in $I(3412, [l_1, \ldots, l_m])$ are symmetric in $l_1,\ldots, l_m$.
\end{enumerate}

%\bibliographystyle{plain}
%\bibliography{references}

\end{document}